\documentclass[11pt]{article}
\usepackage{amsmath,amsthm,amsfonts,amssymb,mathrsfs,epsfig,bm}
\usepackage[usenames]{color}
\usepackage{indentfirst}

\parskip        \smallskipamount

\newtheorem{theorem}{Theorem}%[section]
\newtheorem{lemma}{Lemma}
\newtheorem{corollary}{Corollary}

\newtheorem{note}{Note}
\theoremstyle{remark}

\newtheorem{remark}{Remark}

\def\N{\mathbb{N}}
\def\Z{\mathbb{Z}}

\def\P{\mathbb{P}}
\def\E{\mathbb{E}}
\def\FF{\mathscr{F}}

\def\SS{\mathscr{S}}
\def\PP{\mathscr{P}}
\renewcommand{\phi}{\varphi}
\renewcommand{\epsilon}{\varepsilon}
\newcommand{\1}{{\text{\Large $\mathfrak 1$}}}
\newcommand{\comp}{\raisebox{0.1ex}{\scriptsize $\circ$}}

\newcommand{\var}{\operatorname{var}}

\newcommand{\eqdist}{\stackrel{\text{d}}{=}}

\definecolor{mygray}{gray}{0.9}

\newcommand{\fig}[5]{
\begin{figure}[#1]
\centering
\epsfig{file=#4,width=#2}
\begin{minipage}{#3}
\caption{\em #5}
\end{minipage}
\end{figure}
}
\def\indep{\perp\!\!\!\perp}
\def\before{\bm <\!\!\!\!\;\!\!\bm <}
\def\after{\bm >\!\!\!\!\;\!\!\bm >}

\begin{document}

\title{\bf Limit theorems for a random directed slab graph
\let\thefootnote\relax\footnotetext{Research supported in part by EPSRC grant
EP/E033717/1 and by the Isaac Newton Institute for Mathematical Sciences.
}
\let\thefootnote\relax\footnotetext{
%05C80    	Random graphs
%60F17    	Functional limit theorems; invariance principles
%60K35    	Interacting random processes; statistical mechanics type models; percolation theory
%06A06    	Partial order, general
}
}
\author{\sc Denis Denisov, Serguei Foss \\
\sc and \\
\sc Takis Konstantopoulos\thanks{Corresponding author}}
\date{\small \em \today}
\maketitle

\begin{abstract}
We consider a stochastic directed graph on the integers whereby
a directed edge between $i$ and a larger integer $j$ exists
with probability $p_{j-i}$ depending solely on the distance between
the two integers. Under broad conditions, we identify a regenerative
structure that enables us to prove limit theorems for the maximal path length
in a long chunk of the graph. 
%Alternatively, we speak
%of a stochastic partial order on $\Z$ and a maximal total order.
The model is an extension of a special case of
graphs studied in \cite{FK03}. We then consider a similar type of graph
but on the  `slab' $\Z \times I$, where $I$ is a finite
partially ordered set. We extend the
techniques introduced in the in the first part of the paper
to obtain a central limit theorem for the longest path.
When $I$ is linearly ordered, the limiting distribution can be
seen to be that of the largest eigenvalue of a $|I| \times |I|$
random matrix in the Gaussian unitary ensemble (GUE).
\end{abstract}

\section{Introduction}
Consider a random directed graph with vertex $V=\Z$,
the integers. A pair of integers 
$(i,j)$ is declared to be an edge, directed from $i$ to $j$,
with probability $p_{j-i}$ which depends only on the difference $j-i$,
and this is done independently from pair to pair. 
We assume that $p_{k} =0$ for all $k \le 0$, so there are no directed
edges from a larger integer to a smaller one.
We are interested in limit theorems (law of large number
and central limit theorem) for the maximum length $T[1,n]$ of 
all paths from 1 to $n$, as $n \to \infty$.
The problem as such is related to last-passage percolation.

Unlike nearest-neighbour graphs \cite{MARTIN06,BAI08}, the quantity $T[1,n]$
does not have a direct subadditive property. It turns out that,
a related quantity, namely the maximum
$L[1,n]$ of all paths in the restriction of the graph on $\{1,\ldots,n\}$,
has an almost sub-additive property (see \eqref{L-sub}) 
and thus $L[1,n]/n \to C$, almost surely,
for some deterministic constant $C \le 1$. 
It is later shown that any two vertices are almost surely
eventually connected by a path, and thus $T[1,n]$ has the same asymptotic 
properties as $L[1,n]$.
The minimal condition we need to carry out our programme is
\[
\sum_{k=1}^\infty (1-p_1) \cdots (1-p_k) < \infty.
\]
Under this condition, we can identify a random subset $\SS$ (we call it
``skeleton'') of $\Z$ whose points form 
a stationary renewal process 
(see Sections \ref{renewal} and \ref{independence})
 over which the graph regenerates
and has the property that any element $v$ of $\SS$ is connected
by a path (directed either towards $v$ or away from it)  to any
other vertex in $\Z$.
The quantity $L[1,n]$ becomes additive over the regenerative set $\SS$
enabling us to prove, under the stronger condition
\[
\sum_{k=1}^\infty k (1-p_1) \cdots (1-p_k) < \infty,
\]
a (functional) central limit theorem. The latter condition implies
finiteness of variance of the longest path between two
successive points of $\SS$. To prove the latter assertion, we provide
a rather non-trivial algorithmic construction of the 
last non-positive element of $\SS$. This construction is related
to the so-called coupling-from-the past method for perfect simulation 
\cite{PW97,FOSTWE98} and is the topic of Section \ref{algo}
which is based on the properties of two stopping times
studied in Section \ref{twos}.
The central limit theorem is proved in Section \ref{clt}.

We then consider an extension of the random graph on
the vertex set $\Z \times I$,
where $I$ is a partially
ordered set under some partial order 
$\preceq$ possessing a minimum and a maximum element. 
We let an edge from
$(x,i)$ to $(y,j)$ exist with probability that depends
on $y-x$ and on $i$ and $j$, 
and only when $y-x>0$ and $i \preceq j$. We let $L_N$ be the
length of the longest path in the restriction of the graph on
$\{0,\ldots,N\} \times i$ and show that the
law of $L_N$, appropriately normalized, satisfies a functional central
limit theorem such that the limit process $(Z_t, t \ge 0)$
is $1/2$--self-similar, non-Gaussian, continuous process
with $Z_1$ having the law of the
largest eigenvalue of an a $|I|\times |I|$ random matrix in the Gaussian
Unitary Ensemble (GUE) \cite{AGZ}.

The case where all the $p_k$ are equal to $p$
corresponds to a directed version
of the classical Erd\H{o}s-R\'enyi graph \cite{BE84}. Indeed, let $G_{n,p}$
be the Erd\H{o}s-R\'enyi graph on the set of vertices $\{1, \ldots, n\}$.
To each $\{i,j\}$ which is an edge in $G_{n,p}$ we give an 
orientation from $i \wedge j$ to $i \vee j$. The directed graph thus
obtained is precisely the restriction of our graph on the set 
$\{1, \ldots, n\}$.
This model was also studied in \cite{FK03}. In this paper, we obtained,
among other things, sharp estimates for
the $C \equiv C(p)$ as a function of $p$.
Besides purely mathematical interest,
this model is motivated by applications in Mathematical
Biology (community food webs) \cite{NEWCOH86,COHNEW91,NEWM92},
in Computer Science (parallel processing systems) \cite{ISONEW94},
and in Physics.
Allowing the connectivity probability to depend on the distance between
two vertices $i$ and $j$ means larger modelling flexibility on one hand
while making the model more realistic on the other.

In \cite{FK03} we 
developed a generalisation of Borovkov's theory of renovating events
\cite{BOR78,BOR80,BOR98,BORFOS92,BORFOS94} 
in order to construct a Markov chain in
infinite dimensions describing the ``weights'' of vertices. As a matter of 
fact, in \cite{FK03}, the random graph was a special case of a
more general dynamical system (the ``infinite bin model'') with
stationary and ergodic input. In this paper,
we follow a different approach, one that is applicable specifically for
cases where there is independence between links. In such a case, the approach
has the advantage that it is more elementary using, essentially, renewal theory
and coupling between renewal processes.

\section{The line model}
We are given a set of numbers $(p_j, j \in \N)$, such that 
\[
0 \le p_j < 1, \quad j \in \N.
\]
and consider $(\alpha_{i,j},~ i, j \in \Z,~ i < j)$ as a collection of i.i.d.\
random variables with common law
\[
\P(\alpha_{0,1} =1) = 1-\P(\alpha_{0,1}=-\infty) = p_{j-i}.
\]
Based on this collection, we build a directed random graph $G$
on $\Z$ with edges
\[
E = \{(i,j) \in \Z \times \Z:~ i < j,~ \alpha_{i,j}=1\}.
\]
We shall occasionally refer to the restriction $G[i,j]$ of the graph
on the vertex set $\{i, i+1, \ldots, j\}$ (deleting all edges with either
of the endpoints not in this set).
We are interested in the behaviour of longest paths.
A path $\pi$ is an increasing sequence of vertices
$\pi=(i_0, i_1, \ldots, i_\ell)$ successively connected by edges, i.e.\
$\alpha_{i_0,i_1}=\cdots=\alpha_{i_{\ell-1},i_\ell}=1$. The number
$\ell=|\pi|$ of edges is the length of this path. 

For any $\ell \ge 1$ and any
 increasing sequence $(i_0, i_1, \ldots, i_\ell)$ of vertices
we conveniently define
\begin{equation}
\label{pathlength}
|(i_0, i_1, \ldots, i_\ell)| = 
(\alpha_{i_0,i_1}+ \alpha_{i_1,i_2}+ \cdots
+ \alpha_{i_{\ell-1},i_\ell})^+.
\end{equation}
Clearly, this quantity is $0$ if one of the summands takes value $-\infty$;
otherwise, it equals $\ell$. 
In other words, $|(i_0, i_1, \ldots, i_\ell)| > 0$ if and only if 
$(i_0, i_1, \ldots, i_\ell)$ is a path.

We say that there is a path 
from $i$ to $j$ if $i_0=i$, $i_\ell=j$; we denote this event by
$i \leadsto j$ and may also express it by saying that $i$ leads to $j$
or that $j$ is reachable from $i$.

We let $T[i,j]$ be the maximum length of all paths from $i$ to $j$.
Unlike nearest-neighbour directed graph models 
(see, e.g.\ \cite{MARTIN04}),
this quantity does not have a subadditivity property.
To remedy this we let $L[i,j]$ be the maximum length of all paths from some
$i' \ge i$ to some $j' \le j$, i.e.,
\[
L[i,j] = \max_{i \le i' \le j' \le j} T[i',j'].
\]
That is, $L[i,j]$ is the longest path of the restricted graph $G[i,j]$.
Clearly, $L[i,j]$ has the same law as $L[0,j-i]$.
It is also clear that $L[i,j]$ is subadditive in the sense
that
\begin{equation}
\label{L-sub}
L[i,k] \le L[i,j]+L[j,k] + 1, \quad i < j < k.
\end{equation}
Indeed, if $\pi$ is a path of maximal length in $G[i,k]$ then
its restriction $\pi'$ on $G[i,j]$ has length at most $L[i,j]$ and 
its restriction $\pi''$ on $G[j,k]$ has length at most $L[j,k]$. Now
the length of $\pi$ is equal to the length of $\pi'$ plus the length of
$\pi''$ plus, possibly, 1, if $j$ is not a vertex of $\pi$.
By the subadditive ergodic theorem \cite[p.\ 192]{KALL02},
there exists a deterministic $C \in [0,1]$ such that
\begin{equation}
\label{seg}
\P(\lim_{j \to \infty} L[i,j]/j = C)=1.
\end{equation}

Some of the results below do not depend on the independence assumptions
between the random variables $\alpha_{i,j}$. It is often necessary to
define the model on an appropriate probability space. We do this
as follows.
Let ${\bm \delta} = (\delta_j, j \in \Z)$ be a collection
of independent $\{-\infty,1\}$-valued random variables with
\[
\P(\delta_j=1) =
\begin{cases}
0, & \text{ if } j \le 0
\\
p_j, & \text{ if } j > 0.
\end{cases}
\]
Let ${\bm \delta}^{(i)}$, $i \in \Z$ be i.i.d.\ copies
of $\bm \delta$. The probability space $\Omega$ consists
of $\omega = ({\bm \delta}^{(i)}, i \in \Z)$. 
The random variables $\alpha_{i,j}$ are then defined by
\[
\alpha_{i,j}(\omega) = \delta^{(i)}_{j-i}.
\]
The sigma-field is the standard product sigma-field.
A natural shift $\theta$ on $\Omega$ is the map defined by
\begin{equation}
\label{theta}
%\theta: 
\omega=(i \mapsto {\bm \delta}^{(i)}) \mapsto
\theta\omega=(i \mapsto {\bm \delta}^{(i+1)}).
\end{equation}
Hence
\[
\alpha_{i,j}(\theta \omega) 
= \delta^{(i+1)}_{j-i}
= \delta^{(i+1)}_{(j+1)-(i+1)} = \alpha_{i+1,j+1}(\omega).
\]
The random variables $L[i,j]$ are all defined explicitly on $\Omega$
via $L[i,j] = \max_{i \le i_0 < i_1 < \cdots < i_\ell \le j} 
|(i_0, \ldots, i_\ell)|$ where $(i_0, \ldots, i_\ell)$ is
the random variable defined by \eqref{pathlength}.
It is in this sense that the law $\P$ of the model is $\theta$-invariant
on $\Omega$. Moreover, $\theta$ is ergodic.
In fact, the result that the asymptotic limit of $L[1,n]/n$ exists
depends only on this $\theta$-invariance, so it holds
for more general models where the law of ${\bm \delta}$ is
not that of independent random variables.

A word on notation: If $(A_n, n \in \Z)$ is a collection 
of events of $\Omega$ and $\tau$ is $\Z$-valued random variable on $\Omega$
then $A_\tau$ denotes the event containing all $\omega \in \Omega$
such that $\omega \in A_{\tau(\omega)}$.

\section{The skeleton}\label{renewal}
For the purposes of this section, let $\Omega$ be the space
defined above, $\theta$ the natural shift \eqref{theta}, 
and let $\P$ be a $\theta$-invariant probability measure.
In addition, assume that $\theta$ is ergodic, i.e.\ that 
the invariant sigma-field is trivial.
Recall the shorthand $\{i \leadsto j\} = \{T[i,j]>0\}$
for the event that there is a path from $i$ to $j$. 
%(i.e.\ $j$ is reachable from $i$).
Consider, for each $n \in \Z$, the events 
\begin{align*}
A^+_n&:= \bigcap_{j>n} \{n \leadsto j\}
=  \{\text{any $j>n$ is reachable from $n$}\}\\
A^-_n&:= \bigcap_{j<n} \{j \leadsto n\}
= \{\text{$n$ is reachable from any $j>n$}\}.
\end{align*}
The following is an immediate consequence of the definitions:
\begin{lemma}
(i)
The sequence $\big((A^-_n, A^+_n), n \in \Z\big)$ is stationary
and ergodic.
(ii) For each $n$,
the events $A_n^-$ and $A_n^+$ are independent
and $\P(A_n^+)=\P(A_n^-)=\P(A_0^+)$.
\end{lemma}
We are interested in the random set 
\begin{equation}
\label{skeleton}
\SS(\omega) := \{n \in \Z:~ \omega \in A_n^+ \cap A_n^-\},
\end{equation}
and refer to it as the {\em skeleton} of the random graph.
The terminology is supposed to be reminiscent of a
point of view described next.

%A slightly different point of view then for this model is that
%the random graph $G$ induces a random partial order, namely $\leadsto$,
%on $\Z$. This point of view was taken in \cite{ABBJ94} where linear extensions
%of this random partial order were studied. This point of view is
%important and we shall also keep it in mind below.

Let $\PP(E) \subset \Z\times \Z$ be a partial order (i.e.\ 
if $(i,j), (j,k) \in \PP(E)$ then $(i,k) \in \PP(E)$) which contains
the set of edges $E$. In fact, take $\PP(E)$ to be the smallest such
set. Necessarily, $\PP(E)=\{(i,j) \in \Z \times \Z:~ i \leadsto j\}$.
%Observe that $\leadsto$ is a partial order
%on the set of vertices $V=\Z$. 
%So $i \leadsto j$ means that $i$ precedes
%$j$ in this partial order. 
A subset $U$ of $\Z$ is totally ordered under the partial order $\leadsto$ if
for any distinct $i,j \in U$ we either have $i \leadsto j$ or $j \leadsto i$.
We say that a totally ordered subset $U$ is special if it has the
stronger property that for all distinct $i, j$ with
$i \in U$ and $j \in V$,
we either have $i \leadsto j$ or $j \leadsto i$.
Clearly, the union of special totally ordered subsets is special;
thus we can speak of the maximal special totally ordered subset 
and we refer to it as the skeleton of the partial order.
Adopting this definition,
it is now clear that the set $\SS$ defined by \eqref{skeleton} is
the skeleton of the partial order $\leadsto$ on $\Z$.
In \cite{ABBJ94} the elements of $\SS$ are referred to as {\em posts}.
In fact, \cite{ABBJ94} uses $\SS$ in order to derive limit theorems
of the number $N_n$ of linear extensions of the 
the random partial order $\leadsto$ on $\{1,\ldots,n\}$.

For a general partially ordered set, a skeleton may not exist.
However, in our case, the condition $\P(A_0^+ \cap A_0^-)>0$ is sufficient for
$\SS$ to be almost surely infinite.
\begin{lemma}
If $\lambda:= \P(A_0^+\cap A_0-) >0$ then $\SS$ is an a.s.\ infinite set.
\end{lemma}
\proof
Let $\theta$ be the shift defined by \eqref{theta}. Then,
for all $\omega$, $\SS(\omega) = \SS(\theta\omega)$.
Since $\P$ is $\theta$-invariant, the result follows.
\qed

Assuming that $\lambda = \P(A_0^+\cap A_0^-) >0$, we may then, 
equivalently, consider $\SS$
as a stationary-ergodic point process on the integers with rate
$\lambda$ because $\lambda = \P(0 \in \SS)$.
%\begin{definition}
%For each $\omega \in \Omega$, we say that $n$ is a {\em no-gap vertex}
%if $\omega \in A_n^+ \cap A_n^-$.
%\end{definition}
%The collection of no-gap vertices comprise the random set
%\[
%R(\omega)=\{n \in \Z:~ \omega \in A^+_n\cap A^-_n\}.
%\]
%We now adopt a convention for enumerating the no-gape vertices
%Let
%\[
%\Gamma_1:= \inf(R \cap (0,\infty)),
%\]
%be the first strictly positive no-gap vertex and
%\[
%\Gamma_0:= \sup(R \cap (-\infty, 0)),
%\]
%be the last non-positive no-gap vertex.
%Also, recursively,
%\begin{align*}
%\Gamma_{n+1} &= \min(R\cap \N \setminus \{\Gamma_1, \ldots, \Gamma_n\}),
%\quad n \ge 0
%\\
%\Gamma_{-n} &= \max( R \cap (-\infty,0) 
%\setminus \{\Gamma_1, \ldots, \Gamma_n\}), \quad n \ge 0,
%\end{align*}
%so that 
We let $\Gamma_n$, $n \in \Z$ be an enumeration of the elements  of $\SS$ 
according to the following convention:
\[
\cdots < \Gamma_{-1} < \Gamma_0 \le 0 < \Gamma_1 < \Gamma_2 < \cdots\
\]
In particular, $\Gamma_0$ is the largest non-negative element of $\SS$.
%By stationarity and ergodicity, if $\P(A_0^+)>0$ then
%$\P(|\Gamma_n| < \infty)=1$,
%for all $n \in \Z$.
%By induction, we can define the random sequence 
%\[
%\Gamma_{n+1}=\Gamma_n+\Gamma_1\circ \theta^{\Gamma_n}
%\]
%for $n\in \mathbf Z^+$. Similarly, we define the process $\Gamma_n$ for 
%values of $n\le 0$. 

We can now strengthen the subadditivity property \eqref{L-sub} for $L$:
\begin{lemma}
\label{blocks}
For all integers $m < n$,
$$
L{[\Gamma_m,\Gamma_n]}
=L{[\Gamma_m,\Gamma_{m+1}]}+\cdots+L{[\Gamma_{n-1},\Gamma_n]}.
$$
\end{lemma}
\proof
To see this,
consider the interval $[\Gamma_1,\Gamma_n]$  
and a path $\pi^*$ of length $L{[\Gamma_1,\Gamma_n]}$. 
Then this path must visit all the intermediate
skeleton points $\Gamma_1,\ldots,\Gamma_n$.
Indeed, suppose this is not the case and $\pi^*$ does not visit, say, 
$\Gamma_l$, for some $1\le l\le n$. 
Consider an edge $(i,j)$ belonging to $\pi^*$, with
$i\le \Gamma_l\le j$. 
By the definition of $\Gamma_l$, both
$(i,\Gamma_l)$ and $(\Gamma_l,j)$ are edges of the random graph $G$. 
Therefore we can increase the length of $\pi^*$
by $1$ if we replace the edge $(i,j)$ by two edges 
$(i,\Gamma_l)$ and $(\Gamma_l,j)$. 
This leads to the contradiction since $\pi^*$ has length 
$L{[\Gamma_1,\Gamma_n]}$ which is, by definition, maximal.
\qed

\section{Regenerative structure}\label{independence}
Throughout, we make use of the following two conditions:
\\[5mm]
\hspace*{1cm}
\begin{minipage}{\textwidth}
{\sf [C1] }~ $\displaystyle  0 < p_1 < 1$
\\[5mm]
{\sf [C2] }~ $\displaystyle \sum_{k=1}^\infty (1-p_1)\cdots(1-p_k) < \infty$.
\end{minipage}
\\[5mm]
%\begin{align}
%& 0 < p_1 < 1
%%\label{C1}
%\\
%& \sum_{k=1}^\infty (1-p_1)\cdots(1-p_k) < \infty
%%\label{C2}
%\end{align}
We also sometimes write $q_j = 1-p_j$.
For each $j \in \Z$ we consider its immediate neighbours:
\begin{align}
\overline\eta(j) &:= \min\{k>j:~ \alpha_{j,k} =1\} \nonumber\\
\overline\xi(j) &:= \max\{i<j:~ \alpha_{i,j} =1\}.
\label{nb}
\end{align}
See Figure 1.
The distances of these vertices from $j$ are denoted as follows:
\begin{align*}
\eta(j) &:= \overline\eta(j)  - j \\
\xi(j) &:= j - \overline\xi(j).
\end{align*}
Notice that $(\xi(j), j \in \Z)$ and $(\eta(j), j \in \Z)$ are identically
distributed sequences, and that
each one is a sequence of i.i.d.\ random variables.
Furthermore, for each $j \in \Z$,
\[
(\xi(j+1), \xi(j+2), \ldots)
\indep
(\eta(j-1), \eta(j-2), \ldots)
\]
Henceforth, we shall let $\xi$ be a random variable
with distribution the common distribution of $\xi(j)$ and $\eta(j)$:
\[
\P(\xi > n) 
= \P(\xi(0) > n) 
= \P(\eta(0) > n) 
= (1-p_1) \cdots (1-p_n), \quad n \in \N.
\]
% {\tiny\green
% To make use of this representation we prove two additional facts.
% First we show  that the  $L{[{\Gamma_{k-1}},\Gamma_k]}$ 
% are i.i.d.\  random variables. Second we  show that 
% $\E L{[{\Gamma_{k-1}},\Gamma_k]^2}<\infty$, under appropriate conditions.
% After that, the problem of proving the central limit theorem for  $L_n$ 
% can be solved using Donsker's invariance 
% principle for i.i.d.\ random variables. 
% 
% We are going to realise this program in the subsequent subsections. 
% In Subsection~\ref{independence}  we prove independence of 
% these random variables. In Subsection~\ref{moments} 
% we give conditions for $\P(A^+_0)>0$ and 
% finiteness of the second moment $\E L{({\Gamma_{k-1}},\Gamma_k]^2}$. 
% Finally, in Subsection~\ref{clt} we prove the central limit theorem.
% }

%\comment{
%Comment: We can think of the model in an alternative way:
%Let $P(E) \subset \Z \times \Z$ be the smallest partial order
%containing set of edges $E \subset \Z \times \Z$. 
%Thus, $(i,j) \in P(E)$ means that $i \leadsto j$. It is transitive:
%if $i \leadsto j$ and $j \leadsto k$ then $i \leadsto k$.
%We are looking to identify a subset $R$ of $\Z$ which is totally ordered.
%In fact, the $R$ we defined is precisely it. And it is maximal. It can
%contain no more vertices.
%}

\fig{h}{9cm}{12cm}{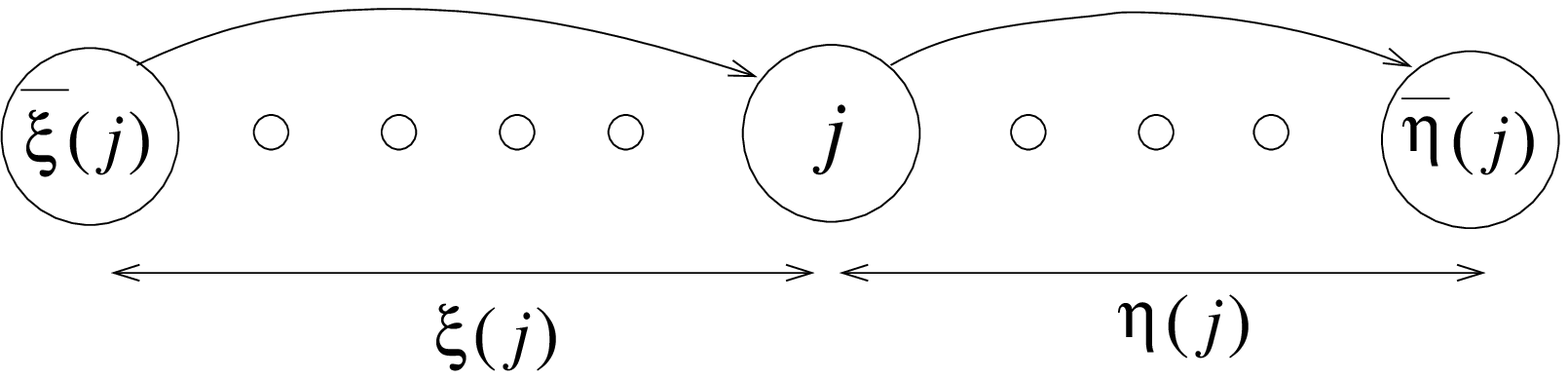}
{Notation used: $\overline \xi(j)$ is the first vertex below $j$ that is
connected to $j$; correspondingly, $\overline \eta(j)$ is the first vertex
above $j$ connected to $j$.}

Define next the events
\begin{equation}
\label{auv}
A_{u,v}^+ := \bigcap_{j=u+1}^v \{ u \leadsto j\}, \quad
A_{u,v}^- := \bigcap_{j=u}^{v-1} \{ j \leadsto v\},
\end{equation}
for which, clearly, 
\[
A^{+}_{u,v} \supset A^{+}_{u,v+1}, \quad
A^{-}_{u,v} \supset A^{-}_{u-1,v}
\]
with
\begin{equation}
\label{AAA}
\lim_{v \to \infty} A_{u,v}^+ = A_u^+, \quad
\lim_{u \to -\infty} A_{u,v}^- = A_v^-.
\end{equation}
Furthermore,
\begin{equation}
\label{BBB}
A_{u,v}^+ \cap A_{v,w}^+ \subset A_{u,w}^+, 
\text{ if } u < v < w,
\end{equation}
a property we shall use in Section \ref{algo}.
Observe also the following:
\begin{lemma}
\label{simple}
For all integers $u<  v$,
\begin{align*}
A_{u,v}^+ & = \bigcap_{j=u+1}^v \bigcup_{i=u}^{j-1} \{i \leadsto j\}
= \bigcap_{j=u+1}^v \{u \le \overline\xi(j)\}
\\
A_{u,v}^-  &= \bigcap_{j=u}^{v-1} \bigcup_{i=j+1}^v \{j \leadsto i\}
 = \bigcap_{j=u}^{v-1} \{\overline\eta(j)\le v\}
\\
A_u^+ &= \bigcap_{j > u} ~\bigcup_{i=u}^{j-1} \{i \leadsto j\}
= \bigcap_{j > u}  \{u \le \overline\xi(j)\}
\\
A_v^- &= \bigcap_{j < v} ~\bigcup_{i=j+1}^{v} \{j \leadsto i\}
= \bigcap_{j < v} \{\overline\eta(j)\le v\}.
\end{align*}
\end{lemma}
\proof
We prove the first equality. %The second follows similarly.
%Recall the definition of $A_n^+ = \bigcap_{j > n} \{n \leadsto j\}$.
%Clearly, $\{n \leadsto j\} \subset \bigcup_{i=n}^{j-1} \{i \leadsto j\}$
%so one inclusion is obvious. 
That $A_{u,v}^+ \subset \cap_{j=u+1}^v \cup_{i=u}^{j-1} \{i \leadsto j\}$
is immediate from the definition \eqref{auv}.
To prove the opposite inclusion, assume that $u > v+1$ (otherwise there
is nothing to prove) and that for all integers
$j \in [u+1, v]$ there exists an integer
$i \in [u, j-1]$ such that $i \leadsto j$. Fix $j > u $ and pick $i_1$ to be
the largest among the vertices between $u$ and $j-1$ such that $i_1 \leadsto
j$; necessarily, $\alpha_{i_1,j}=1$. Then pick the largest vertex $i_2$
among the vertices between $u$ and $i_1-1$ such that $i_2 \leadsto i_1$,
and continue this way.
Since $i_1 > i_2 > \cdots \ge u$, it follows that this process terminates
with some $i_k =u$. Since $(u=i_k, i_{k+1}, \ldots, i_1, j)$ is
a path, we have that $u \leadsto j$.
The second equality for $A_{u,v}^+$ now follows from the definition
\eqref{nb}.
The relations for $A_{u,v}^-$ follow similarly.
The third (respectively, fourth) line is obtained by sending 
$v$ to $+\infty$ (respectively, $u$ to $-\infty$)
in the first (respectively, second) one. 
\qed

This lemma tells us that $A_{u,v}^+$ is the intersection of
$v-u$ independent events. Indeed, since $\overline\xi(j) = j-\xi(j)$
we have
\begin{equation}
\label{Auv+}
A_{u,v}^+ =\{\xi(u+1) \le 1, \xi(u+2) \le 2, \ldots, \xi(v) \le v-u\}, 
\end{equation}
and the random variables $\xi(u+1),\ldots,\xi(v)$ are i.i.d.
Similarly, for $A_{u,v}^-$, 
\begin{equation}
\label{Auv-}
A_{u,v}^- = \{\eta(u) \le v-u, \ldots, \eta(v-2) \le 2, \eta(v-1)\le 1\}. 
\end{equation}
Moreover, since
\[
(\xi(u+1), \xi(u+2), \ldots, \xi(v)) \eqdist
(\eta(v-1), \eta(v-2), \ldots, \eta(u))
\]
we have that $\P(A_{u,v}^+) = \P(A_{u,v}^-)$.
Similarly, both $A_n^+$ and $A_n^-$ are intersections 
of infinitely many independent events:
\begin{align}
A_n^+ 
&= \bigcap_{j>n} \{\xi(j) \le j-n\}
\label{An+}
\\
A_n^- 
&= \bigcap_{j<n} \{\eta(j) \le n-j\},
\label{An-}
\end{align}
and $\P(A_n^+)= \P(A_n^-)$.
The skeleton \eqref{skeleton} can be expressed as follows:
\begin{equation}
\label{Sexpr}
\SS %= A_n^- \cap A_n^+ 
= \{n \in \Z:~
%\ldots, \overline\eta(n-2), \overline\eta(n-1) 
\sup_{i < n} \overline\eta(i)
\le 
n
\le
\inf_{j > n} \overline\xi(j)\}.
%\overline\xi(n+1), \overline\xi(n+2),\ldots\}
\end{equation}
Regarding $\SS$ as a point process, we see that it has rate
\[
\lambda = \P(0 \in \SS) = \P(A_0^+)^2 
%= \P( \cap_{j \ge 1} \{\xi_j \le j \})
= \big( \prod_{j=1}^{\infty} \P(\xi(j)\le j) \big)^2
= \prod_{j=1}^{\infty} [1- \P(\xi(0) > j)]^2.
\]
Since
\begin{equation}
\label{xilaw}
\P(\xi(0) > j) = \P(\alpha_{0,1} = \cdots = \alpha_{0,j}=0)
= (1-p_1) \cdots (1-p_j),
\end{equation}
we have
\begin{equation}
\label{lambdaformula}
\lambda = \prod_{j=1}^{\infty} 
[1- (1-p_1) \cdots (1-p_j)]^2
\end{equation}
and so
\[
\text{{\sf [C2]} } \iff \lambda>0 \iff
\E[\xi(0)] < \infty.
\]

%Consider now the infinite-dimensional vectors
%\[
%\bm \alpha_n := (\alpha_{n,n+1},~ \alpha_{n,n+2}, \ldots),
%\quad n \in \Z.
%\] 
%We shall refer to them as {\em connectivity vectors}
%and show that the process $(\bm \alpha_n, n \in \Z)$ regenerates
%on the no-gap vertices. Define the cycles of the process over
%successive no-gap vertices by
\newcommand{\alphacycle}{\mathscr{C}}
Consider now two successive skeleton points $\Gamma_k$ and $\Gamma_{k+1}$
and let $\alphacycle_k(\omega)$ be the restriction of $\omega$ on $[\Gamma_k,
\Gamma_{k+1})$:
\[
\alphacycle_k
:= \big(\bm \delta^{(n)},~ \Gamma_k \le n < \Gamma_{k+1}\big),
\quad k \in \Z;
\]
we refer to it as the $k$-th ``cycle''.
We next show that the sequence of cycles 
have a regenerative structure in the following sense:
\begin{lemma}
\label{reg}
The cycles $(\alphacycle_k, k \in \Z)$ are independent and
$(\alphacycle_k, k \in \Z-\{0\})$ are identically distributed.
In particular, the skeleton vertices $(\Gamma_k, k \in \Z)$ 
form a stationary renewal process.
\end{lemma}
Intuitively, Lemma \ref{reg} is based on the following observation.
Suppose that $0$ is a skeleton vertex (i.e.\ condition on the event
$A_0^- \cap A_0^+$). Then $\overline\xi(1) \ge 0$,
$\overline\xi(2) \ge 0$, etc.
In other words, $\overline\xi(1)=0$, $\overline\xi(2)\in \{1,2\}$,
$\overline\xi(3) \in \{0,1,2\}$, etc. To determine the location of the
next skeleton vertex after $0$ we need to find the first vertex $j > 0$
such which is connected with every vertex between $0$ and $j-1$.
This means that, conditional on $0$ being a skeleton vertex,
the location of the first skeleton vertex larger than $0$
does not depend on the $(\bm \delta^{(n)}, n < 0)$.

\noindent
{\em Proof of Lemma \ref{reg}.}
\\
For $k \ge 1$, 
let $\FF_k^+$ be the sigma-algebra generated by $(\bm\delta^{(1)}, \ldots,
\bm\delta^{(k)})$ and let $\FF_k^-$ be the sigma-algebra
generated by $(\bm\delta^{(-1)}, \ldots, \bm\delta^{(-k)})$. 
It suffices to prove that, for any $k \ge -1$, $l \ge 1$, and any
$B_k^- \in \FF_k^-$, $B_l^+ \in \FF_l^+$,
\begin{equation*}
\P (\Gamma_{-1}=k,B^-_k,\Gamma_1=l,B^+_l \mid \Gamma_0=0)
=\P (\Gamma_{-1}=k,B^-_k \mid \Gamma_0=0)
\P (\Gamma_1=l,B^+_l \mid \Gamma_0=0).
\end{equation*}
Assume that $\Gamma_0=0$ (i.e.\ $0$ is a skeleton vertex).
Then, by \eqref{Sexpr},
%\begin{align*}
%\Gamma_1 &= \min \{ n>0  :~  \1_{A_n^-\cap A_n^+}=1 \},\\
%\Gamma_{-1} &= \max \{ n<0  :~  \1_{A_n^+\cap A_n^-}=1 \}.
%\end{align*}
\[
\ldots,~ \overline\eta(-2),~ \overline\eta(-1)
\le 0 \le \overline\xi(1),~ \overline\xi(2),~ \ldots
\]
In view of the latter inequality, we have
\begin{align*}
\Gamma_{-1} &= \max \{ n<0  :~  \1_{A_n^-\cap A_n^+}=1 \}
\\ 
&= \max \{ n < 0:~ \ldots,~ \overline\eta(n-2),~ \overline\eta(n-1)
\le n \le \overline\xi(n+1),~ \overline\xi(n+2),~ \ldots\}
\\
&= \max \{ n < 0:~ \ldots,~ \overline\eta(n-2),~ \overline\eta(n-1)
\le n \le \overline\xi(n+1),~ \overline\xi(n+2),~ \ldots,~ \overline\xi(0)\}
=: \widehat \Gamma_{-1},
\end{align*}
where the last serves as a definition of a new random variable 
$\widehat\Gamma_{-1}$.
This random variable
is $\FF^-$-measurable, where $\FF^-$ is the sigma-algebra
generated by $(\bm\delta^{(k)}, k <0)$.
Similarly, we define
\[
\widehat\Gamma_1:=\min\{n > 0:~ \overline\eta(0), \ldots, \overline\eta(n-1)
\le n \le \overline\xi(n+1), \overline\xi(n+2), \ldots\},
\]
a random variable which is $\FF^+$-measurable, 
where $\FF^+$ is the sigma-algebra
generated by $(\bm\delta^{(k)}, k >0)$,
and observe that, on $\{\Gamma_0=0\}$,
the random variables $\Gamma_1$ and $\widehat\Gamma_1$ coincide.
%\begin{align*}
%\Gamma_1 &= \min \{ n>0  :~  \1_{A_n^-\cap A_n^+}=1 \}
%\\
%&= \min\{n > 0:~ \overline\eta(0), \ldots, \overline\eta(n-1)
%\le n \le \overline\xi(n+1), \overline\xi(n+2), \ldots\},
%\end{align*}
%and the latter is $\FF^+$-measurable, where $\FF^+$ is the sigma-algebra
%generated by $(\bm\alpha_k, k >0)$.
Note that $\FF^-$ and  $\FF^+$ are independent. 
Hence,  for $k \le -1$, $\ell \ge 1$, we have
\begin{align*}
\P (\Gamma_{-1}=k,B^-_k,\Gamma_1=l,B^+_l &| \Gamma_0=0)
\\
&= \frac{\P (\Gamma_{-1}=k,B^-_k,\Gamma_1=l,B^+_l, A_0^+\cap A_0^-)}
{\P (A_0^+\cap A_0^-)}\\
&= \frac{\P (\{\widehat\Gamma_{-1}=k\}\cap A_0^-\cap B^-_k , 
\{\widehat\Gamma_1=l\}\cap A_0^+\cap B^+_l)}{
\P (A_0^+)\P(A_0^-)}\\
&= \frac{\P(\{\widehat\Gamma_{-1}=k\}\cap A_0^-\cap B^-_k) 
~\P( \{\widehat\Gamma_1=l\}\cap A_0^+\cap B^+_l)}{
\P (A_0^+)\P (A_0^-)}
\\
&= \P(\widehat\Gamma_{-1}=k,  B^-_k \ | \ A_0^-)
~\P( \widehat\Gamma_1=l,B^+_l \ | \ A_0^+).
\end{align*}
Note that
\begin{multline*}
\P(\widehat\Gamma_{-1}=k ,  B^-_k \ | \ A_0^-)
= \P(\widehat\Gamma_{-1}=k,  B^-_k \ | \ A_0^-\cap A_0^+)
\\
= \P(\Gamma_{-1}=k,  B^-_k \ | \ A_0^-\cap A_0^+)
=\P(\Gamma_{-1}=k,  B^-_k \ | \Gamma_0=0).
\end{multline*}
Similarly,
$$
\P \{ \widehat\Gamma_{1}=1 ,  B^+_l \ | \ A_0^+ \}
=\P \{\Gamma_{1}=1,  B^+_l \ | \Gamma_0=0\}.
$$
%Therefore 
%conditioned on the event $\{\Gamma_0=0\}$ the random vectors 
%$(\gamma_0,\chi_0)$ and $(\gamma_1,\chi_1)$ are i.i.d.
%The induction arguments (which we omit to avoid cumbersome notation) imply the general statement that $\{(\gamma_l,\chi_l)\}_{l=-\infty}^\infty$ is a sequence of i.i.d. vectors.
\qed
\begin{corollary}
\label{iidstruktur}
The bivariate random variables 
\[
\big(\Gamma_1-\Gamma_0,~ L[\Gamma_0, \Gamma_1]\big),~~
\big(\Gamma_2-\Gamma_1,~ L[\Gamma_1, \Gamma_2]\big),~~ \ldots
\]
are i.i.d.
\end{corollary}

\section{Two stopping times}
\label{twos}
In this section, we study properties of the following two random variables:
\begin{align*}
\mu &:= \inf\{i >0:~ \1_{A_{-i,0}^-}=0\}\\
\nu &:= \inf\{i >0:~ \1_{A_{-i,0}^+}=1\}.
\end{align*}
These random variables are important in the algorithmic
construction of Section \ref{algo}.

Note that $-\nu$ is the first vertex $< 0$ 
with the property that every vertex in the interval $(-\nu, 0]$ is
reachable from $-\nu$:
\[
\nu = \inf\{i>0:~ -\nu \leadsto 0,~ -\nu \leadsto -1,~ \ldots,~
-\nu \leadsto -\nu+1\}.
\]
Also, $-\mu$ is the first vertex $< 0$ such that 
$0$ is not reachable from $-\mu$:
\[
\mu = \inf\{ i > 0:~ -i \not \leadsto 0\}.
\]

We will show that $\mu$ is a defective random variable,
i.e.\ that $\P(\mu = \infty) > 0$, with conditional tail
$\P(\mu > n | \mu < \infty)$ comparable to the
integrated tail of $\xi$. We will also show that $\nu$ is
an a.s.\ finite random variable with the same number of moments
as $\xi$.

Note first that both $\mu$ and $\nu$ are stopping times with respect to the
filtration $(\FF_{k}^-, k \le 0)$.
Observe that
\begin{equation}
\label{muinfty}
\{\mu=\infty\} = \bigcap_{i \ge 1} A_{-i,0}^- = A_0^-.
\end{equation}
Since condition {\sf [C2]} is equivalent to $\P(A_0^-)>0$, we have
\[
\P(\mu=\infty) > 0.
\]
On the other hand,
\[
\{\nu = \infty\} = \bigcap_{n=1}^\infty (A_{-n,0}^+)^c,
\]
and, as we shall see below, this event has probability zero:
\begin{equation}
\label{nuinfty}
\P(\nu = \infty) = 0.
\end{equation}

% {\green
% Next recall that 
% \[
% A_0^- = \{ \sup_{j \ge 1} \overline\eta(-j) \le 0\},
% %= \{ \eta(-j) \le j \text{ for all } j \ge 1\},
% \]
% and let
% \[
% \mu := \inf\{j \ge 1:~ \overline \eta(-j) >0\}
% = \inf \{j\ge 1 :~ \ \eta(-j)>j \}
% \]
% Then
% \[
% A_0^- = \{ \mu=\infty \}.
% \]
% }

Let us first focus on the law of $\mu$, conditional on $\{\mu < \infty\}$.
This can be computed easily, from the definition of $\mu$, and equations
\eqref{Auv-}, \eqref{muinfty}, and \eqref{xilaw}.
\begin{align}
\P(n < \mu < \infty)
&= \P(\eta(-k) \le k \text{ for all } 1\le k \le n)~
\P(\eta(-m) > m \text{ for some } m > n) 
\label{muexact}
\\
&= \prod_{k=1}^n \P(\eta(-k) \le k)~
\bigg(1-\prod_{m=n+1}^\infty \P(\eta(-m) \le m) \bigg)
\nonumber
\\
&=(1-q_1)(1-q_1q_2) \cdots (1-q_1q_2\cdots q_n)
\bigg(1-\prod_{m=n+1}^\infty (1-q_1q_2\cdots q_m) \bigg)
\nonumber
\end{align}
Conditional on $\{\mu < \infty\}$, the random variable $\mu$
has a tail comparable to the integrated tail of $\xi$:
\begin{lemma}\label{lem1}
Suppose that {\sf [C1]} and {\sf [C2]} hold. 
There exist constants $0<C_1<C_2<\infty$ such that, for all $n \ge 0$,
\[
C_1 \sum_{m>n}^{\infty} \P (\xi >m)
\le 
\P (\mu > n \ | \ \mu < \infty)
\le 
C_2 \sum_{m>n}^{\infty} \P (\xi >m).
\]
\end{lemma}
\proof
Since $p_1 < 1$, we have $\lambda < 1$ (see \eqref{lambdaformula}) and so
\[
\P(\mu <\infty )= 1-\lambda^{1/2} > 0.
\]
%The upper bound is obtained as follows: from the the definition of $\mu$
%and \eqref{Auv-}, \eqref{An-} we have
%\[
%\{n < \mu < \infty\}
%= \{\eta(-1) \le 1, \ldots, \eta(-n) \le n\} \cap
%\{\eta(-m) > m \text{ for some } m > n\}.
%\]
%so that
%\[
%\P(n < \mu < \infty)
%= (1-q_1)(1-q_1q_2)\cdots(1-q_1q_2\cdots q_n)
%\sum_{m>n} \P(\xi > m).
%\]
Using \eqref{muexact} we have
\begin{align*}
\P(\mu > n \ | \ \mu < \infty)
%&= \frac{1}{1-\lambda^{1/2}}~ \P(n<\mu <\infty) 
%\\
%&= \frac{1}{1-\lambda^{1/2}}~
%\P(\eta(m)\le m, 1\le m\le n) \cdot
%\P(\eta(m)>m \text{ for some } m>n)
%\\
&\le  \frac{1}{1-\lambda^{1/2}}~\sum_{m=n+1}^\infty\P(\eta(m)>m)
=  \frac{1}{1-\lambda^{1/2}}~\sum_{m=n+1}^\infty\P(\xi>m).
\end{align*}  
Hence $C_1 = 1/(1-\lambda^{1/2})$.
To obtain a bound from below note that 
the first term on the right of $\eqref{muexact}$ is $\ge \P(\mu=\infty)$
and so
\begin{align*}
\P(n < \mu < \infty) 
%&= \P(\eta(-k) \le k \text{ for all } 1\le k \le n)~
%\P(\eta(-m) > m \text{ for some } m > n)
%\\
%&\ge \P(\eta(-k) \le k \text{ for all }  k \le n)~
%\P(\eta(-m) > m \text{ for some } m > n)
%\\
&= \lambda^{1/2}
\bigg(1-\prod_{m=n+1}^\infty \P(\eta(-m) \le m) \bigg)
\\
& \ge
\lambda^{1/2}
\bigg(1-\exp \big(-\sum_{m=n+1}^\infty \P(\xi>m) \big) \bigg)
\\
& \ge
\lambda^{1/2} g(\E\xi) \sum_{m=n+1}^\infty \P(\xi>m),
\end{align*}
where $g(x) = (1-e^{-x})/x$.
Hence $C_2 = g(\E\xi) \lambda^{1/2}/(1-\lambda^{1/2}))$.
\qed

We next prove something stronger than \eqref{nuinfty}, namely that
$\nu$ has the same number of moments as $\xi$.
\begin{lemma}\label{lem2}
If $\E\xi^r < \infty$ for some $r \ge 1$ then $\E\nu^r <\infty$.
\end{lemma}
\proof
By the definition of $\nu$ and equation \eqref{Auv+} we have
\[
\nu =\inf\{n\ge 1:~\xi(0)\le n,\xi(-1)\le n-1,\ldots,\xi(-(n-1))\le 1\}
\]
Define a sequence of non-negative random variables $x_0, x_1, x_2, \ldots$ by
$x_0=0$ and
\[
x_n=\max\{\xi(0)-n, \xi(-1)-(n-1),\ldots, \xi(-(n-1))-1\}, \quad n \ge 1.
\]
Then
\[
\nu=\inf\{ n\ge 1:~ x_n=0 \}.
\]
The $x_n$ satisfy 
\[
x_{n+1}=\max(x_n,\xi(-n)) -1, \quad n \ge 0,
\]
and, since the $\xi(-n)$ are i.i.d., 
$(x_n, n \ge 0)$ is a Markov chain in $\Z_+$.
We now make two observations that imply the statement of the lemma. 
First, if  $x_n>K>0$ then 
\[
x_{n+1}-x_n
= (\xi(-n)-x_n)^+ -1
\le (\xi(-n)-K)^+ -1.
\]
But $\E[(\xi-K)^+]<1$ for sufficiently large $K$. 
Therefore, after the Markov chain leaves the interval $[0,K]$
(for sufficiently large $K$) 
it is majorized from above by a random walk with increments distributed like
$(\xi-K)^+-1$ whose mean is negative. 
By standard properties of random walks this implies that 
the return time $T_K$ to the set $[0,K]$ satisfies $\E T_K^r < \infty$
if $\E((\xi-K)^+-1)^r<\infty$; and the latter is equivalent to 
$\E \xi^r < \infty$.
The second observation is that the Markov chain 
$(x_n)$ returning to the set $[0,K]$  
eventually hits point $0$ after a geometric number of trials.  
\qed

\begin{corollary}
\label{ccc}
If {\sf [C2]} holds then $\E \nu < \infty$.
\end{corollary}

%\comment{
%The rest needs to be polished and fused with the previous.
%There is repetition... But, for now, I will not worry about this.
%Some notation which might be useful is as follows.
%Let $u < v$ be integers.
%\begin{align*}
%A^+_{u,v} &:= \{ u \leadsto j ~ \forall j=u+1, \ldots, v\}
%%\{u \text{ leads to every vertex in } [u,v]\}
%\\
%A^-_{u,v} &:= \{ j \leadsto v ~ \forall j=u, \ldots, v-1\}
%= \{\text{every vertex in } [u,v] \text{ leads to } v\}
%\end{align*}
%Equivalently,
%\begin{align*}
%A^+_{u,v} &:= \{ u \le \overline \xi(u+1), \ldots, \overline \xi(v) \}
%= \{ \xi(u+1) \le 1, \ldots, \xi(v) \le v-u \}
%\\
%A^-_{u,v} &:= \{ \overline\eta(u), \ldots, \overline\eta(v-1) \le v\}
%= \{ \eta(u) \le v-u , \ldots, \eta(v-1) \le 1\}
%\end{align*}
%}

\section{Algorithmic construction of $\Gamma_0$}
\label{algo}
In this section we give a method for constructing a specific skeleton point,
e.g., the first one which is to the left of the origin. This is the point 
$\Gamma_0$. Besides the theoretical interest, such a construction will be
used later for proving a central limit theorem; it can also be used
in connection to a perfect simulation algorithm for estimating the
value of $C=\lim_{n \to \infty}L[1,n]/n$ (see remarks at the end of the section).

The idea for the construction of $\Gamma_0$ is this: recall that
$-\nu$ which is the first vertex 
$< 0$ which is connected to every point between $-\nu$ and $0$.
We check whether $-\nu$ is also reachable from every point from the left.
If it is, we declare that $-\nu$ is a {\em silver} point and stop the procedure.
If not, there is a first vertex before $-\nu$ which fails to be 
connected to $-\nu$.
Using the shift operator $\theta$ defined in \eqref{theta}, this
vertex is at distance $\mu\comp \theta^{-\nu}$ from $-\nu$;
in other words, this distance is the functional $\mu$ applied to the shifted
$\omega$, when the origin is placed at $-\nu$.
We then set $\mu[1] = \nu+\mu \comp \theta^{-\nu}$, which is the location
of the previous vertex, and $\nu[1]=\nu$ and this finishes the first
step of the procedure.

The second step of the algorithm is similar to the first one:
we search for the first vertex $-\nu[2]$ before $-\mu[1]$ which is connected
to every vertex between $-\nu[2]$ and $-\nu[1]$. We know that we can find
such a vertex with probability one. If it also happens that $-\nu[2]$
is reachable from any point from the left, we stop and declare
$-\nu[2]$ as our {\em silver} point. Otherwise, there will be a first
vertex, $-\mu[2] < -\nu[2]$ which fails to be connected to $-\nu[2]$.

The procedure continues in the same way, until the first {\em silver} point 
is found, and it will be found with probability one.
This first silver point will have the property that it is reachable from every
point from the left and is connected to every point up until the origin; see
Lemma \ref{halflemma} below.
The distribution of this first silver point
is well-understood and this is the content of Lemma \ref{silverdist}.
In fact, we will show that there are infinitely many silver points which form
a (delayed) renewal process backwards; see Lemma \ref{silverseq}.
Finally, in Theorem \ref{gamma0} we show that among the infinitude of
silver points we can pick a {\em gold} one, namely the point $\Gamma_0$.

To define the algorithm explicitly,
we consider a sequence of $\N\cup\{+\infty\}$-valued stopping times relative to
the filtration 
$(\FF_k, k \ge 1)$, defined as follows.
Let
\begin{align}
\nu[1] &:=\nu
\nonumber\\
\mu[1] &:= \nu +  \mu \comp \theta^{-\nu}
= \inf\{j > \nu:~ \1_{A^-_{-j,-\nu} }=0 \},
\label{r1}
\end{align}
and, recursively, for $k \ge 2$,
\begin{align}
\nu[k] &:= 
\inf\{j > \mu[k-1]: \1_{A_{-j, -\nu[k-1] }^+  }=1 \}
\nonumber\\
\mu[k] &:= \nu[k] + \mu \comp \theta^{-\nu[k]}
= \inf\{j > \nu[k]: \1_{A_{-j,-\nu[k]}^-  } =0\},
\label{r2}
\end{align}
where $\theta$ is the natural shift \eqref{theta}.
It is understood that if for some $k$
we have $\mu[k]=\infty$ then $\nu[j]=\mu[j]=\infty$ for all $j \ge k+1$.
We thus obtain an increasing sequence of stopping times
\[
\nu=\nu[1] < \mu[1] < \nu[2] < \mu[2] < \nu[3] < \mu[3] < \cdots
\]
which (since $\P(\mu=\infty)>0$) is eventually equal to infinity.
It is convenient to think of these stopping times as the points
of an alternating point process (the $\mu$-points and the $\nu$-points).
In words, the sequence of these stopping times is defined by
first laying a $\nu$-point in location $\nu[1]$. Then, as long
as $\eta(-(\nu[1]+i)) \le i$ for $i=1,2,\ldots$, we place
no point in location $\nu[1]+i$. At the first instance $i$ at which
$\eta(-(\nu[1]+i)) > i$, we place a $\mu$-point
in location $\nu[1]+i$ and call it $\mu[1]$. 
The random variables $(\eta(-(\nu[1]+i)), i \ge 1)$ 
are independent of $\nu[1]$, and so the event that we place a $\mu$-point
in a finite location is independent of
$\nu[1]$ and has probability $\P(\mu<\infty) = 1-\lambda^{1/2}$.
The procedure continues in the same way: having placed $\nu[k]<\infty$,
we decide, independently of the past (i.e.\ $\FF_{\nu[k]}^-$) whether
to create a new $\mu$-point or not (i.e.\ place it at infinity).
If we do create a new $\mu$-point $\mu[k]$ then, clearly, $\nu[k+1]$
is also finite and $\nu[k+1]-\nu[k]$ has the 
same distribution as $\nu[2]-\nu[1]$ conditional on $\mu[1]<\infty$.
Thus for each $\omega$, the recursion stops at the index
\begin{equation}
\label{K}
K:= \inf\{k\ge 1:~ \mu[k]=\infty\}.
\end{equation}
{From} the discussion above we immediately obtain:
\begin{lemma} 
Assume that {\sf [C1]} and {\sf [C2]} hold. Then $K$ is a geometric
random variable with
\[
\P(K >k) = (1-\lambda^{1/2})^k, \quad k \ge 0.
\]
\end{lemma}
By definition, $\mu[K]=\infty$ but $\mu[K-1]<\infty$. Hence
\[
\nu[K] < \infty, \text{ a.s.}
\]

\begin{note}
We stop for a minute to point out that the whole purpose
of the construction of these random variables is the random variable
$\nu[K]$. In other words, for each $\omega \in \Omega$, we apply
recursion \eqref{r1}-\eqref{r2} to obtain the alternating sequence
of $\nu$ and $\mu$- points, through them we define that index
$K$ as in \eqref{K} and, finally, $\nu[K]$. Thus, $\nu[K]$ is
a well-defined (measurable) function of $\omega$.
We refer to $-\nu[K]$ as the first {\em silver} point before $0$.
\end{note}

Although $K$ depends on the whole alternating process 
$(\nu[k],\mu[k]), k \ge 1)$, we can identify the law of $\nu[K]$
as follows:
\begin{lemma}
\label{silverdist}
On a new probability space, let $K, \psi_1, \psi_2, \psi_3, \ldots$
be independent random variables with distributions
\begin{align*}
& P(K > k) = (1-\lambda^{1/2})^k, \quad k \ge 0
\\
&\psi_1 \eqdist \nu 
\\
& \psi_{i}  \eqdist \big(\nu[2]-\nu[1] \mid \mu[1]< \infty)
\eqdist \big(\inf\{j > \mu :~ \1_{A_{-j,0}^+}=1\} \mid \mu < \infty\big),
\quad i \ge 2.
\end{align*}
Then, assuming {\sf [C1]} and {\sf [C2]},
\begin{equation}
\nu[K] \eqdist \psi_1 + \sum_{i=1}^{K-1} \psi_{i+1}.
\label{nudist}
\end{equation}
\end{lemma}
\proof
It follows from
\[
\nu[K] = \nu[1] + \sum_{i=1}^{K-1} (\nu[i+1]-\nu[i]).
\]
using a simple probabilistic argument as described above.
\qed

%{\tiny\green
%Perhaps we need a complete proof here.}

The reason we are interested in the random variable $\nu[K]$ is the following:
\begin{lemma}
\label{halflemma}
Assume {\sf [C1]} and {\sf [C2]} hold. Then for $\P$-a.e.\ $\omega$
\begin{equation}
\label{halfconnected}
\omega \in A_{-\nu[K]}^- \cap A_{-\nu[K],0}^+.
\end{equation}
\end{lemma}

Note that replacing the index $n$ in a sequence of events
$A_n$ by a random index $N$ amounts to defining the event
$A_N=\{\omega \in \Omega:
\text{ there exists } n \text{ such that } n = N(\omega) \text{ and }
\omega \in A_n \}$.

The meaning of \eqref{halfconnected} is that the vertex $\nu[K]$
of the random graph has the property that there is a path from
every $j < \nu[K]$ to $\nu[K]$ and there is a path
from $\nu[K]$ to every $i$ such that $\nu[K] < i \le 0$.
Our goal is to identify a skeleton point. Whereas $\nu[K]$ is not
a skeleton point for sure, there is a positive probability that it is.

\noindent
{\em Proof of Lemma \ref{halflemma}.}
%\begin{align*}
%\nu[k] &= \mu[k-1] + \inf\{i > 0: \1_{A_{-i,0}^+}=1\} \comp \theta^{-\mu[k-1]}
%\\
%&= \mu[k-1] + \inf\{i > 0:  \1_{A_{-(i+\mu[k-1]),-\mu[k-1]}^+}=1\}
%\\
%&= \inf\{j > \mu[k-1]: \1_{A_{-j,-\mu[k-1]}^+}=1\}
%\end{align*}
%\begin{align*}
%\{K < \infty\} 
%&\subset \{\mu[K]=\infty\} \cap \{\mu[K-1]<\infty, \ldots, \mu[0] < \infty\}
%\\
%&\subset \{\mu[K]=\infty\} \cap \{\nu[K]<\infty, \ldots, \nu[0] < \infty\},
%\end{align*}
%by Corollary \ref{ccc}.
If $K=k$, for some $k \ge 1$,  then
$\mu[k]=\infty$ but $\mu[k-1] < \infty$, so
$\nu[k] < \infty$ and
$\1_{A_{-j,-\nu[k]}^-}=0$ for all $j > \nu[k]$.
Hence
\[
%\{\mu[k]=\infty, \mu[k-1] < \infty\}
\{K=k\}
\subset \bigcap_{j > \nu[k]} A_{-j,-\nu[k]}^- 
= A_{-\nu[k]}^-,
\]
by \eqref{AAA}.
Also, if $K=k$, then $\nu[k], \nu[k-1], \ldots, \nu[1] < \infty$
and so
\[
\{K=k\} \subset A^+_{-\nu[k],-\nu[k-1]} \cap 
A^+_{-\nu[k-1],-\nu[k-2]} \cap \cdots \cap A^+_{-\nu[1],0}
\subset A^+_{-\nu[k],0},
\]
by \eqref{BBB}.
But $K$ is a geometric random variable 
and hence $K<\infty$, a.s.
\qed

We also have the following result concerning moments of $\nu[K]$:
\begin{lemma}
Assume {\sf [C1]} and {\sf [C2]} hold.
If, in addition, there exists
$r \ge 1$ such that $\E \xi^{r+1}<\infty$, then $\E\nu[K]^r < \infty$.
\end{lemma}
\proof
We have that $\E\nu[K]^r < \infty$ if $\E \nu^r < \infty$
and $\E(\mu^r | \mu < \infty) < \infty$.
The latter holds if $\E \xi^{r+1}<\infty$, and this is
a simple consequence of Lemma \ref{lem1}. On the other hand,
$\E \nu^r < \infty$ holds if $\E \xi^r < \infty$, as proved in
Lemma \ref{lem2}.
\qed

Whereas {\sf [C1]} and {\sf [C2]} imply $\P(\nu[K]<\infty)$,
we need finite variance for $\xi$ in order that we have
finite expectation for $\nu[K]$.

% {\tiny\green
% In this section we provide a recursive construction of the 
% the set $R \cap (-\infty, 0]$ and, as a corollary, obtain
% further conditions for the moments of the typical distance
% between two no-gap vertices, i.e.\ the two vertices straddling
% the origin.
% 
% \comment{
% Throughout, we assume that 
% \[
% 0< 1-r := \P(A_0^-) < 1.
% \]
% The second inequality fails only in a trivial case. 
% For the first inequality we gave conditions.
% }
% 
% The construction is achieved in two steps.
% 
% {\em Step 1:} We wish to identify vertices $v <0$ such that
% $\1_{A_v^- \cap A_{v,0}^+}=1$, i.e.\ vertices $v$ such that there is a path 
% from every vertex $u< v$ to $v$ and there is a path from $v$ to every
% $v<w\le 0$. 
% 
% {\em Step 2:} Amongst the vertices $v$ found by Step 1, we select
% those for which $\1_{A_v^+}=1$. (Note that $A_v^+ \subset A_{v,0}^+$.)
% These comprise the set $R \cap (-\infty, 0]$.
% 
% 
% ====================% 
% STUFF DELETED
% 
% =====================% }

We next construct a further sequence of stopping times.
\[
\sigma[1] < \sigma[2] < \cdots
\]
as follows.
Assume that {\sf [C1]} and {\sf [C2]} hold.
Recall that the random variable $\nu[K]$ 
is a.s.\ finite; it maps $\Omega$ into $\N$.
Hence we can define $\nu[K]\comp\theta^{n}$ for any $n \in \Z$ and
also $\nu[K]\comp \theta^{J}$ for any measurable $J: \Omega \to \Z$.
We define $\sigma[j]$, $j \ge 1$, recursively:
\begin{align}
& \sigma[1] = \nu[K] 
\nonumber
\\
& \sigma[j+1] = \sigma[j] + \nu[K] \comp \theta^{-\sigma[j]},
\quad j \ge 1.
\label{sigma}
\end{align}
Intuitively, given $\omega$,
we first construct $\nu[K]$ by \eqref{r1}-\eqref{r2} and place a point 
$\sigma[1]$ at $\nu[K]$.  
We then shift the origin to $-\nu[K]$ and repeat the recursion 
with $\omega' = \theta^{-\nu[K]}(\omega)$ in place
\footnote{$\omega' = \theta^{-\nu[K(\omega)](\omega)}(\omega)$} of $\omega$,
thus obtaining a new random variable, $\nu[K] \comp \theta^{-\nu[K]}$.
We place another point $\sigma[2]$ at distance
\footnote{$\nu[K] \comp \theta^{-\nu[K]}(\omega)
= \nu[K(\omega')](\omega') = \nu[K(\theta^{-\nu[K(\omega)](\omega)}(\omega))](\theta^{-\nu[K(\omega)](\omega)}(\omega))   $}
 $\nu[K] \comp \theta^{-\nu[K]}$
from $\sigma[1]$. The procedure continues in the same way.
We refer to $-\sigma[1], -\sigma[2], \ldots$
as the sequence of silver points.

\begin{lemma}
\label{silverseq}
Assume that {\sf [C1]} and {\sf [C2]} hold.
Define the point process with points $\sigma[j]$, $j \ge 1$, as in
\eqref{sigma}. This is a renewal process on $\N$, i.e.\ the
random variables $\sigma[1]$, $\sigma[2]-\sigma[1]$, $\sigma[3]-\sigma[2], 
\ldots$ are i.i.d.\ with common distribution \eqref{nudist}.
\end{lemma}

We are now ready to construct the first gold point $\Gamma_0$. 
\begin{theorem}
\label{gamma0}
Assume that {\sf [C1]} and {\sf [C2]} hold.
Define the sequence $(\nu[k], \mu[k], k \ge 1)$ through \eqref{r1}-\eqref{r2}
which is used to define the random variable $\nu[K]$.
Based on this, define the sequence $(\sigma[j], j \ge 1)$, through 
\eqref{sigma}. In addition,
let
\begin{align*}
M &:= \sup_{i \ge 1} \{\xi(i)-i\},
\\ 
J &:= \inf\{j \ge 1:~ \sigma[j] \ge M \}.
\end{align*}
Then
\[
\Gamma_0 = -\sigma[J].
\]
\end{theorem}
Before proving the theorem, let us observe that the random
variables defined in the theorem statement are a.s.-finite.
By {\sf [C2]}, i.e.\ that $\E \xi < \infty$, 
implies $M < \infty$, a.s.
\begin{align}
\P(M \ge m) &= \P(\xi(i)-i \ge m, \text{ for some } i \ge 1)
\nonumber
\\
&\le \sum_{i=1}^\infty \P(\xi(i) \ge i+m)
\\
&\le \sum_{i=m+1}^\infty \P(\xi(i) \ge i) \le \E \xi.
\label{Mtail}
\end{align}
By standard renewal theory, it is easy to see that $J$, the first
exceedance of $M$ by the random walk $(\sigma[j], j \ge 1)$, 
is also a.s.-finite and hence $\sigma[J]$ is an a.s.-finite random
variable.

\noindent
{\em Proof of Theorem \ref{gamma0}.}
%The assumptions imply that $\E\xi <\infty$, and this implies that 
%$\E M < \infty$. Hence $\sigma[J] < \infty$, a.s.
Owing to Lemma \ref{halflemma}, we have that
\begin{equation}
\label{sprop}
\text{for all } j \in \N, \quad
\omega \in A_{-\sigma[j]}^- \cap A_{-\sigma[j],0}^+, \quad \P-{a.e.\ }~ \omega
\in \Omega.
\end{equation}
Also, 
\begin{equation}
\label{also}
\{M \le \sigma[J]\} = \{ \xi(1) \le \sigma[J]+1, \xi(2) \le \sigma[J]+2,
\ldots \}.
\end{equation}
Fix $n \in \N$ and observe that, from the definition of $M$ and 
the expressions \eqref{Auv+}, \eqref{An+} for $A^+_{-n,0}$ and $A^+_{-n}$,
respectively,
\begin{align*}
A^-_{-n} \cap A^+_{-n,0} \cap \{ M \le n\}
&= A^-_{-n} \cap A^+_{-n,0} \cap \{ \xi(1) \le 1, \xi(2) \le 2, \cdots \}
\\
&= A^-_{-n} \cap  \{\xi(-n+1) \le 1, \ldots, \xi(0) \le n, \xi(1) \le 1, \xi(2) \le 2, \cdots \}
\\
&= A^-_{-n} \cap   A^+_{-n} 
\\
& = \{ n \in \SS\}.
\end{align*}
Combining this with \eqref{sprop} and \eqref{also} we obtain
\[
-\sigma[J] \in \SS, \quad \text{a.s.}
\] 
It is clear, from the algorithmic construction \eqref{r1}-\eqref{r2} of
the sequence $(\nu[k], \mu[k], k \ge 1)$, from the algorithmic construction
\eqref{sigma} of the $(\sigma[j], j \ge 1)$, and the definition of $J$,
that there can be no point of $\SS$ between $-\sigma[J]$ and $0$.
Therefore $-\sigma[J]$ is the largest negative point of $\SS$. 
\qed

\begin{remark}
Possible extensions:
The algorithmic construction proposed above may be used in a general stationary
ergodic framework.
In particular, one can easily generalise first-order results 
(the functional strong law of large numbers). 
Under reasonable assumptions, one
can again prove the finiteness of 
$\xi (0)$. This  will imply the finiteness
of $\eta (0)$ and, in turn, the existence of the stationary skeleton. 
Then the functional strong law of large numbers
will follow using well-known tools.
\end{remark}

\begin{remark}
Simulation and perfect (exact) simulation of the value of the limit $C$: 
This depends in a complex way on an infinite number of variables, and 
one cannot expect an analytic closed form expression. 
But one can estimate it by 
running a MCMC algorithm. One can also use the 
regenerative structure of the model 
to run the simulation in backward
time using the idea of ``cycle-truncation'' 
that leads to a simple implementation
scheme; c.f.\ \cite{FTC} for more details
However, each such an algorithm gives a biased
estimator of the unknown parameter, in general.

In \cite{FK03}, we considered the homogeneous case ($p_j=p$, for all $j$).
In particular, in \cite[\S 10]{FK03}
(see also \cite[\S 4]{FK03} for theoretical background), 
we obtained a stronger result by proposing an algorithm
for the perfect simulation of a random sample from an 
unknown distribution whose mean is the limit $C$ under consideration. 
The standard MCMC scheme provides an unbiased
estimator for this limit. 

The ideas behind that algorithm may be efficiently implemented in a number
of similar models, e.g.\ in models with long memory (
see, for example, \cite{COMFERFER}).  
In fact, in \cite{FK03}, we developed the algorithm for a more general model 
(we called it ``infinite-bin model'')
and under general stochastic ergodic assumptions. 
%The latter model is equivalent 
%to our currect model only if the i.i.d.\ assumptions are in force.
\end{remark}

\section{Central limit theorem for the maximum length}\label{clt}
Assume now that {\sf [C1]} holds and
\\[5mm]
\hspace*{1cm}
\begin{minipage}{\textwidth}
{\sf [C3] }~ $\displaystyle \sum_{k=1}^\infty k (1-p_1)\cdots(1-p_k) < \infty$.
\end{minipage}
{From} \eqref{xilaw} we see that this is equivalent to 
\\[5mm]
\hspace*{1cm}
\begin{minipage}{\textwidth}
{\sf [C3$'$] }~ $\displaystyle \E \xi^2 < \infty$. 
\end{minipage}
\begin{lemma}
\label{finitevar}
If {\sf [C1]} and  {\sf [C3]} hold then $\E |\Gamma_0| < \infty$.
\end{lemma}
\proof
By Theorem  \ref{gamma0}, $|\Gamma_0| = \sigma[J]
= \min\{\sigma[j]:~ j \ge 1,~ \sigma[j] \ge M\}$.
Recall that $\sigma[1] < \sigma[2] < \cdots$ are points
of a renewal process. This renewal process is clearly independent 
of $M = \sup_{i \ge 1} \{\xi(i)-i\}$. By standard renewal theory, 
$\E \sigma[J] <\infty$ if $\E M < \infty$. But the tail of $M$ was
estimated in \eqref{Mtail}. The same inequalities now show that
$\E \xi^2 < \infty$ is sufficient for $\E M < \infty$.
\qed

The maximum length $L_n$ of all paths from some $i \ge 0$
to some $j \le n$ satisfies the following central limit theorem.
\begin{theorem}
\label{Lclt}
Suppose {\sf [C1]} and {\sf [C3]} hold.
Let 
\[
\sigma^2 := \var\big(L(\Gamma_1, \Gamma_2] - C (\Gamma_2-\Gamma_1)\big).
\]
Define
\[
\ell_n(t) := \frac{L_{[nt]} - Cnt}{\lambda^{1/2}\sigma \sqrt{n}},
\quad t \ge 0, \quad n \in \N.
\]
Then the sequence of processes $\ell_n$, in the Skorokhod
space $D[0,\infty)$ equipped with the topology of uniform
convergence on compacta \cite{BILL68}, converges weakly to a standard
Brownian motion.
\end{theorem}
\proof
By Lemma \ref{finitevar} we have $\E|\Gamma_0|<\infty$.
Hence $\E\Gamma_1 < \infty$. But the $\Gamma_n$ form a stationary
renewal process. Therefore, $\E\Gamma_1 < \infty$ implies that
the variance of $\Gamma_2-\Gamma_1$ is finite.
Since $L(\Gamma_1,\Gamma_2] \le \Gamma_2-\Gamma_1$, 
we have $\sigma^2 < \infty$.
The constant $C$, defined as the a.s.-limit of $L_n/n$--see
\eqref{seg}, is also finite and nonzero.
Lemma \ref{blocks} shows that $(L_n, n \ge 0)$ is a (stationary)
regenerative process. The result then is then obtained by reducing
it to Donsker's theorem. This is standard, but we sketch 
the reduction here for completeness. 
Let $\Phi_n$ be the cardinality of $\SS \cap [0,n]$ (the number of $\Gamma_j$
in the interval $[0,n]$):
\[
\Phi_n:= |\SS \cap [0,n]| = \sum_{j \in \Z} \1(0 \le \Gamma_j \le n).
\]
So $\Gamma_{\Phi_n} \le n < \Gamma_{\Phi_n+1}$.
Write
\begin{align*}
L_{[nt]} &= \{ L_{[nt]} - L_{\Gamma_{\Phi_{[nt]}}} \} + L_{\Gamma_{\Phi_{[nt]}}}
\\
nt &= \{ nt - \Gamma_{\Phi_{[nt]}} \} + \Gamma_{\Phi_{[nt]}}.
\end{align*}
The quantities in brackets on both lines are tight and so they
are negligible when divided by $\sqrt{n}$. 
So instead of $\ell_n(t)$, we consider
\begin{equation}
\label{adj}
\widehat \ell_n(t) := 
\frac{L_{\Gamma_{\Phi_{[nt]}}} - C \Gamma_{\Phi_{[nt]}}}{\lambda^{1/2}\sigma \sqrt{n}}
= 
\frac{L_{\Gamma_1}-C\Gamma_1}{\lambda^{1/2}\sigma \sqrt{n}}
+ \frac{1}{\lambda^{1/2}\sigma \sqrt{n}} 
 \sum_{i=2}^{\Phi_{[nt]}} \{L(\Gamma_{i-1},\Gamma_i]-C(\Gamma_i-\Gamma_{i-1})\}
\end{equation}
The last term is the one responsible for the weak limit of $\widehat \ell_n$
(and hence of $\ell_n$). To save some space, put
\[
\chi_i := L(\Gamma_{i-1},\Gamma_i]-C(\Gamma_i-\Gamma_{i-1}).
\]
Donsker's theorem says that
\[
\bigg(\frac{1}{\sigma\sqrt{n}} \sum_{i=2}^{nu} \chi_i, u \ge 0 \bigg)
\Rightarrow (B_u, u \ge 0),
\]
weakly in $D[0,\infty)$, as $n \to \infty$, where $B$ is a standard
Brownian motion.
Let
\[
\phi_n(t) := \frac{\Phi_{[nt]}}{n}, \quad t \ge 0.
\]
Since $\phi_n$ converges
weakly, as $n \to \infty$, to the deterministic function 
$\big(\lambda t, t \ge 0)$ and since composition
is a continuous operation, the continuous mapping theorem
tells us that 
\[
\bigg(\frac{1}{\sigma\sqrt{n}} \sum_{i=2}^{n \phi_n(t)} \chi_i, u \ge 0 \bigg)
\Rightarrow (B_{\lambda u}, u \ge 0) \eqdist \lambda^{1/2} B,
\]
and this readily implies that the last term in \eqref{adj} converges
weakly to a Brownian motion.
\qed

It is now easy to see how the quantity $T[i,j]$, the maximum length
of all paths from $i$ to $j$, behaves. A sufficient condition
for $T[i,j]$ to be positive is that there is a skeleton point between
$i$ and $j$. Therefore, keeping $i$ fixed, the probability that eventually
for all $j$ sufficiently large $T[i,j] >0$ is at least
equal to the probability that eventually there is a skeleton point in $[i,j]$,
and this is certainly equal to one.
So, eventually, any two points are connected, a.s.

Moreover, 
\[
T[\Gamma_i, \Gamma_j]= L[\Gamma_i, \Gamma_j].
\]
Indeed, $\Gamma_i$ is connected to every larger vertex and
any vertex smaller than $\Gamma_j$ is connected to $\Gamma_j$.
Thus, if a path from some $u \ge \Gamma_i$ to some $v \le \Gamma_j$
has length $L[\Gamma_i, \Gamma_j]$ we necessarily have $u=\Gamma_i$
and $v=\Gamma_j$ and this shows the equality of the last display.

If $n$ is large enough so that there is at least one skeleton
point in $[0,n]$, we have that $0 \leadsto n$ and
\[
L[\Gamma_1, \Gamma_{\Phi_n}] \le T[0,n] \le L[\Gamma_0, \Gamma_{\Phi_n+1}],
\]
where $\Phi_n$ is the number of skeleton points in $[0,n]$.
Therefore we immediately obtain:
\begin{theorem}
If {\sf [C1]} and {\sf [C2]} hold then $T[0,n]/n \to C$, as
$n \to \infty$, a.s.
\end{theorem}

Same rationale shows:
\begin{theorem}
Suppose {\sf [C1]} and {\sf [C3]} hold. Then Theorem \ref{Lclt} holds
with $T$ in place of $L$.
\end{theorem}

\section{Directed slab graph} 
Recall that we started with vertex set $V=\Z$ and introduced a random partial order
$\leadsto$ by means of a random directed graph: 
\begin{equation}
\label{redef}
i \leadsto j \text{ if } i< j \text{ and
$\exists$ } i=i_0 < i_1 < \cdots < i_\ell=j \text{ such that }
\alpha_{i_0,i_1}= \cdots = \alpha_{i_{\ell-1},j}=1.
\end{equation}
A natural generalisation is to replace the
total order $<$ on the vertex set $V$
by a partial order $\prec$ and substitute the $i< j$ requirement in
\eqref{redef} above by
the requirement that $i \prec j$. 
We here provide an example of such a generalisation.
A major role in our analysis has been played by the
assumption that the underlying probability measure is invariant by some
shift $\theta$.  Our example will also satisfy this assumption.

Let $(I, \preceq)$ be a finite partially ordered set.
We assume that $I$ has a minimum and a maximum, denoted by
$0$ and $M$, respectively. In other words, for all $i, j, k \in I$,
\begin{align*}
(a)\quad & 0 \preceq i \preceq i \preceq M,
\\
(b)\quad&  \text{if } i \preceq j \preceq i \text{ then } i=j,
\\
(c)\quad & \text{if } i \preceq j \preceq k \text{ then } i \preceq k.
\end{align*}
Consider $V = \Z \times I$. 
We call this vertex set a {\em cylinder}. In the case
$I = \{0,1,\ldots, M\}$, with the usual ordering, we call $V$ a {\em slab}.
Elements of $V$ will be denoted by $(x,i)$, $(y,j)$, etc.
We introduce the component-wise partial ordering $\before$ on $V$ by
\[
(x,i) \before (y,j) \iff (x,i) \neq (y,j) \text{ and } 
x \le y,~  i \preceq j,
\]
and write $(y,j) \after (x,i)$ for the same thing.
Next, we assign an edge $\big((x,i),(y,j)\big)$ to each
pair of vertices such that $(x,i) \before (y,j)$ with probability
$r_{y-x, i,j}$, independently from pair to pair.
This is done by means of random variables $\alpha_{(x,i),(y,j)}$:
\[
\P(\alpha_{(x,i),(y,j)} = 0) = 1-\P(\alpha_{(x,i),(y,j)} =-\infty)r_{y-x,i,j}.
\]
We shall make this more formal in the sequel.
The problem is, again, the behaviour of a longest path from $(x,i)$
to $(y,j)$. This length is denoted by $T[(x,i),(y,j)]$.
We also define $L[(x,i),(y,j)]$ to be the maximum
length of all paths starting from some $(x',i') \after (x,i)$ and ending
at some $(y',j') \after (y,j)$.

An appropriate probability space for the model is now described. Let
$\bm \delta = (\delta_{x,i,j}, x \in \Z, i,j \in I)$ be a collection of 
independent $\{-\infty,1\}$-valued random variables with
\[
\P(\delta_{x,i,j}=1)=r_{x,i,j},
\]
assuming that $r_{x,i,j}=0$ if $x \le 0$ or if $i \succ j$.
Next, let $\bm \delta^{(x)}$, $x \in \Z$ be a collection of i.i.d.\ copies
of $\bm \delta$. The probability space $\Omega$ is defined to contain
infinite vectors $\omega=(\bm \delta^{(x)}, x \in \Z)$.
In other words, $\Omega = (\{-\infty,1\}^{\Z \times I \times I})^\Z$
with $\{-\infty,1\}^{\Z \times I \times I}$ be the space of values
of each $\bm \delta^{(x)}$, and with $\P$ being a product measure.
A shift $\theta$ on $\Omega$ is taken to be the natural
map
\begin{equation}
\omega=(x \mapsto {\bm \delta}^{(x)}) \mapsto
\theta\omega=(x \mapsto {\bm \delta}^{(x+1)}).
\end{equation}
Clearly, $\P$ is preserved by $\theta$.
The random variables $\alpha_{(x,i),(y,j)}$ are now given by
\[
\alpha_{(x,i),(y,j)}(\omega) = \delta^{(x)}_{y-x,i,j}
\]
and it is easy to check their $\theta$-compatibility:
$\alpha_{(x,i),(y,j)}(\theta \omega) = \alpha_{(x+1,i),(y+1,j)}(\omega)$.

We introduce the following assumptions on the probabilities $r_{x,i,j}$.
\\[5mm]
\hspace*{1cm}
\begin{minipage}{\textwidth}
{\sf [D0] }~ 
$\displaystyle r_{x,i,i} =: p_x$ for all $i \in I$
\\[5mm]
{\sf [D1] }~ 
$\displaystyle  0 < p_{1} < 1$
\\[5mm]
{\sf [D2] }~ 
$\displaystyle \sum_{x=1}^\infty  (1-p_{1})\cdots(1-p_{x}) 
< \infty$
\\[5mm]
{\sf [D2$'$] }~ 
$\displaystyle \sum_{x=1}^\infty x  (1-p_{1})\cdots(1-p_{x}) 
< \infty$
\\[5mm]
{\sf [D3] }~ 
For all $i,j \in I$ with $i \prec j$, we have $r_{0,i,j} > 0$.
\end{minipage}

Of these, the last one is not an essential condition. It is only
introduced for convenience. We will comment on it later.
Of course, {\sf [D2$'$] } is stronger than {\sf [D2] } and
it will be used for the proof of the CLT.

\subsection{The random graph $G[x,y]$}
The random directed graph $G=(V,E)$ with $V= \Z \times I$ and $E$
consisting of all $\big((x,i),(y,j)\big)$ such that $\alpha_{(x,i),(y,j)}=1$
is now a well-defined object. 
Let $G[x,y]$ be the restriction of $G$ on the vertex set
$[x,y] \times I$ where $x \le y$ are two integers.
Let
\[
L[x,y] := \max_{\substack{x \le x' \le y' \le y\\i,j \in I}} L[(x',i),(y',j)]
\]
be the maximum length of all paths in $G[x,y]$.
We have $\theta$-compatibility
\[
L[x,y] \comp \theta = L[x+1,y+1],
\]
and, by an argument analogous to the one used to obtain \eqref{L-sub},
we have the subadditivity property
\[
L[x,z] \le L[x,y] + L[y,z]+1, \quad x \le y \le z.
\]
Therefore,
\[
L_N/N := L[0,N]/N \to C, \text{ as $n \to \infty$, a.s.},
\]
for some deterministic constant $C$ which, under the assumption {\sf [D2]},
is positive.

\subsection{The random graph $G^{(i)}$}
Let $G^{(i)}$ be the restriction
of $G$ on the vertex set $V \times \{i\}$, $i \in I$.
It is clear that each $G^{(i)}$ is a line model as studied earlier.
In fact, the $G^{(i)}$, $i \in I$ are i.i.d.
We denote by $L^{(i)}[x,y]$ the maximum length of all paths
of $G^{(i)}$ from some vertex $x' \ge x$ to some vertex $y' \le y$.
We shall let $\SS^{(i)}$ be the skeleton of $G^{(i)}$.
Then, assuming {\sf [D1]} and {\sf [D2]},
each $\SS^{(i)}$  forms a stationary renewal process with nontrivial rate.
Moreover, {\sf [D1]} implies that this renewal process is aperiodic.

\section{Central limit theorem for the directed cylinder graph}
We first describe the limiting process.
To do this, we need the following.
First, let $(B^{(i)}(t), t \ge 0)$, $i \in I$, be i.i.d.\ standard Brownian
motions, all starting from $0$.
Second, let $H(I,\preceq)$ be the {\em Hasse diagram} \cite{DP90}
corresponding 
to the partially ordered set $I$. This is a directed graph
with vertex set $I$ and an edge from $i$ to $j$ if there
is no $k$, distinct from $i$ and $j$, such that $i \preceq k \preceq j$.
Let $\iota=(\iota_0, \iota_1, \ldots, \iota_r)$ be a path
in $H(I,\preceq)$ starting from $\iota_0=0$ and ending at $\iota_r=M$.
The length of the path is $r=|\iota|$.
For each such path $\iota$, define the stochastic process 
$(Z(\iota)_t, t \ge 0)$ by:
\begin{equation}
\label{Z1}
Z(\iota)_t :\sup_{\substack{0 \le t_0 \le t_1 \le \cdots \le t_{|\iota|}=t}}
\big\{
B^{(\iota_0)}(t_{0})
+ [ B^{(\iota_1)}(t_{1})-B^{(\iota_1)}(t_{0}) ]
+ \cdots
+ [ B^{(\iota_{|\iota|})}(t_{|\iota|})-B^{(\iota_{|\iota|})}(t_{|\iota|-1}) ]
\big\}
\end{equation}
and then let
\begin{equation}
\label{Z2}
Z_t := \max_\iota B(\iota)_t,
\end{equation}
where the maximum is taken over all paths $\iota$ 
from the minimum to the maximum element in the Hasse diagram.

The main theorem of this section is as follows:
\begin{theorem}
\label{ccclt}
Let $G$ be a directed cylinder graph and assume that 
{\sf [D0]}, {\sf [D1]}, %{\sf [D2]}, 
{\sf [D2$'$]}, {\sf [D3]}
hold. Let $L_n$ be the maximum length of all paths in $G[0,n]$.
There exists a constant $\kappa > 0$ such that 
\[
\ell_n(t) := \frac{L_{[nt]} - Cnt}{\kappa \sqrt{n}},
\quad t \ge 0, \quad n \in \N
\]
converges weakly, as $n \to \infty$, in the Skorokhod
space $D[0,\infty)$ equipped with the topology of uniform
convergence on compacta, to the stochastic process $Z$ defined in 
\eqref{Z1}-\eqref{Z2}.
\end{theorem}
\proof
Since the $\SS^{(i)}$, $i \in I$ are independent aperiodic renewal processes,
we have 
that
\[
\SS := \{x \in \Z:~ x \in \cap_{i\in I} \SS^{(i)}, ~
\alpha_{(x,i),(x,j)}=1 \text{ for all $i,j \in I$ with $i \prec j$}\}
\]
is also a renewal process. 
Indeed, Lindvall \cite{LIN92} shows that $\cap_i \SS^{(i)}$ is
a stationary renewal process. Now $\SS$ is obtained
from $\cap_{i\in I} \SS^{(i)}$ by a further independent
thinning with positive probability due to the convenient assumption {\sf [D3]}.
Condition {\sf [D2]} implies that
the rate of each $\SS^{(i)}$ is positive and this implies that the
rate of $\cap_{i\in I} \SS^{(i)}$ is positive. 
Hence the rate of $\SS$ is also positive.
Call this rate $\lambda$. We have $0< \lambda \le  1$. Moreover,
$\SS$ is stationary: $\SS \comp\theta = \SS$.
Enumerate now the points of $\SS$ by
\[
\cdots <  \Gamma_{-1} < \Gamma_0 <0 \le \Gamma_1 < \Gamma_2 < \cdots
\]
\fig{b}{9cm}{12cm}{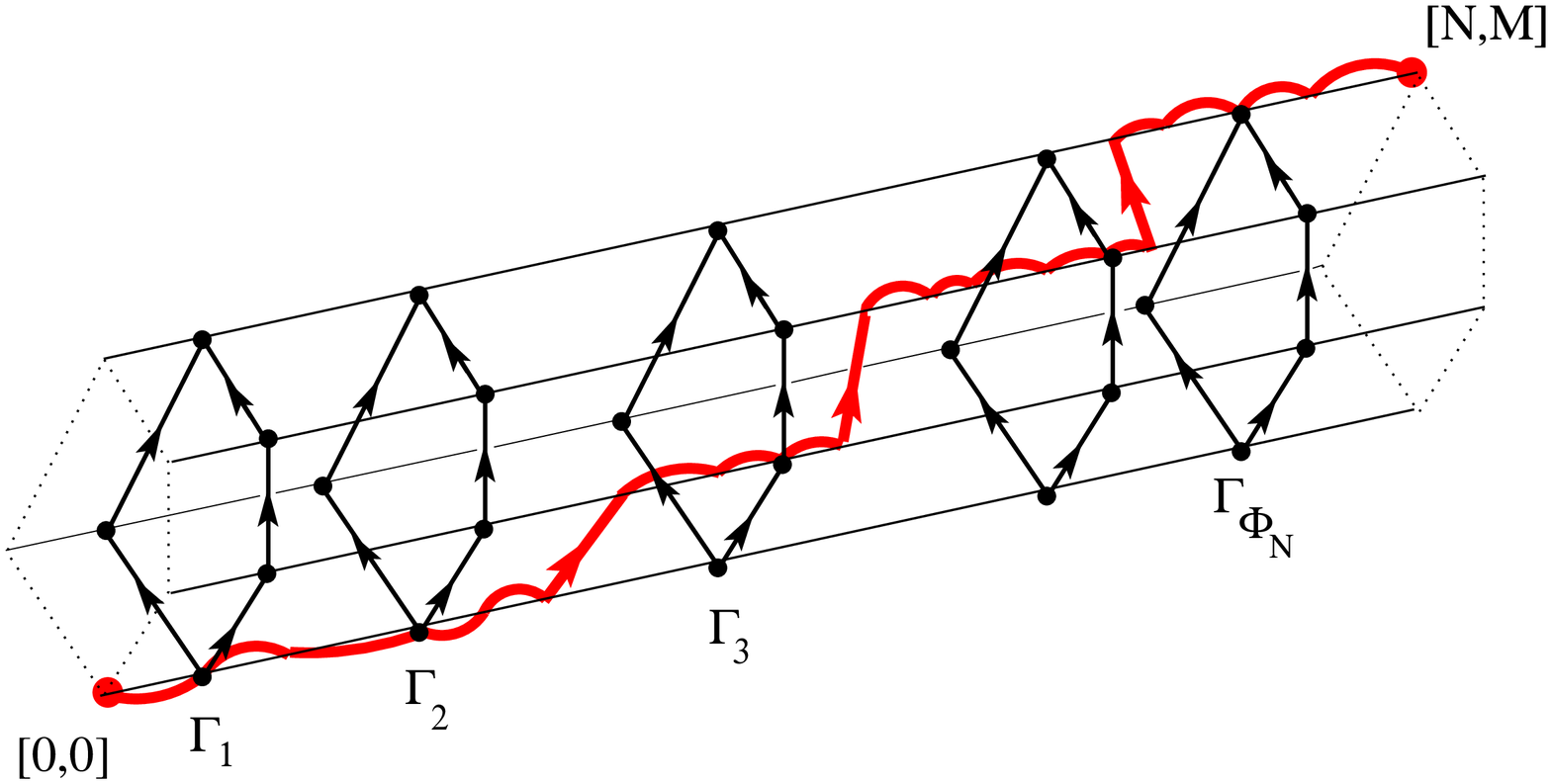}
{The skeleton for the slab graph and a longest path.}
We have $\E (\Gamma_2-\Gamma_1) = 1/\lambda$. If 
\[
\Phi_n := | \SS \cap [0,n]|,
\]
we have $\lim_{n \to \infty}\Phi_n/n =\lambda$, a.s.
Furthermore, $C=\lambda \E L[\Gamma_1, \Gamma_2] \le 1$.
Condition {\sf [D2$'$]} implies that $\E(\Gamma_2-\Gamma_1)^2<\infty$
and hence $\E L^{(i)}[\Gamma_2-\Gamma_1]^2 < \infty$.
By Corollary \ref{iidstruktur}, the random variables
\[
\big(\Gamma_2-\Gamma_1,~ L^{(i)}[\Gamma_1, \Gamma_2]\big),~~ 
\big(\Gamma_3-\Gamma_2,~ L^{(i)}[\Gamma_3, \Gamma_2]\big),~~
\ldots
\]
are i.i.d., and since $\SS$ is obtained by
independent thinning of $\cap_{i\in I} \SS^{(i)}$,
we further have that the rows of the last display are also independent
when $i$ ranges in $I$.

Consider next a path
$\iota=(\iota_0, \iota_1, \ldots, \iota_r)$,  of length $|\iota|=r$
in the Hasse diagram $H(I, \preceq)$ and define the quantities
\begin{align*}
L^*(\iota)_n &:= \max_{1\le j_0 \le j_1 \le \cdots \le j_r =n}
\{
L^{(\iota_0)}[\Gamma_1, \Gamma_{j_0}] +
L^{(\iota_1)}[\Gamma_{j_0}, \Gamma_{j_1}] +
\cdots
+ 
L^{(\iota_r)}[\Gamma_{j_{r-1}}, \Gamma_{j_r}]
\}
\\
L^*_n &:= \max_\iota L^*(\iota)_n,
\end{align*}
where the last maximum is taken over all paths $\iota$ from the minimum
to the maximum element of the Hasse diagram.

We now argue that the quantity of interest $L_n$ is of order 
$L^*_n + o_{\mathrm  d}(1)$
when $n$ is large by providing an upper and a lower bound.
The key observation is that when $n$ is large, the number of points
$\Gamma_j \le n$ grows at a positive rate (and hence to infinity).
At each of these points, say $\Gamma_j$, the graph $G[\Gamma_j, \Gamma_j]$
(being a vertical slice of $G$--see Figure 2) 
is precisely the Hasse diagram:
\[
G[\Gamma_j, \Gamma_j] = H(I, \preceq), \quad j \in \Z.
\]
Fix $\iota' \prec \iota''$ in $I$.
Since $\Gamma_j$ is a point in the skeleton of $G^{(\iota')}$,
any $x \le \Gamma_j$ is connected to $\Gamma_j$ in $G^{(\iota')}$.
Similarly, $\Gamma_j$ is connected to any $y$ in $G^{(\iota'')}$.
Since $\iota'$ is connected to $\iota''$ in $G[\Gamma_j, \Gamma_j]$,
it follows that, almost surely, there is path in $G$ from any $(x,\iota')$ to
any $(y,\iota'')$, if $x \le \Gamma_j \le y$ for some $\Gamma_j \in \SS$
and if $\iota' \prec \iota''$.

Assume that $\Phi_n \ge 2$. Let $\iota=(\iota_0, \iota_1, \ldots, \iota_r)$
be a path in $H(I, \preceq)$ with $\iota_0=0$, $\iota_r=M$ and consider
integers 
\begin{equation}
\label{jjj}
1 \le j_0 \le j_1 \le \cdots \le j_{r-1} \le j_r=\Phi_n.
\end{equation}
Keep in mind that 
\[
\Gamma_{\Phi_n} \le n.
\]
By the construction of the set $\SS$, the following is true:
\[
(\Gamma_1,0) = (\Gamma_1, \iota_0) \leadsto (\Gamma_{j_0}, \iota_0)
\leadsto (\Gamma_{j_0}, \iota_1) \leadsto (\Gamma_{j_1}, \iota_1) \leadsto
\cdots \leadsto (\Gamma_{j_{r-1}},\iota_r) \leadsto (\Gamma_{j_r}, \iota_r)
= (\Gamma_{\Phi_n}, M),
\]
where $(x,\iota') \leadsto (y,\iota'')$ means that there is a path from
$(x,\iota')$ to $(y,\iota'')$ in $G$.
Therefore
\[
L_n \ge L^{(\iota_0)}[\Gamma_1, \Gamma_{j_0}] +
L^{(\iota_1)}[\Gamma_{j_0}, \Gamma_{j_1}] +
\cdots
+
L^{(\iota_r)}[\Gamma_{j_{r-1}}, \Gamma_{j_r}],
\]
because the right-hand side is a lower bound on the
length of the specific path chosen in the last display.
By keeping $\iota$ fixed and maximising over the $j_0, \ldots, j_r$
satisfying \eqref{jjj}
we obtain $L_n \ge L^*(\iota)_n$, and by maximising over $\iota$
we obtain the lower bound 
\[
L_n \ge L^*_{\Phi_n}.
\]

To obtain an upper bound, let $\pi^*$ be a path that achieves
the maximum in $L_n$. Assume that $\Phi_n \ge 1$ so that,
by the key observation above,
$(0,0)$ is connected to $(n,M)$ in $G$. See Figure 3.
Hence $\pi^*$ is necessarily a path from $(0,0)$ to $(n,M)$.

%The elements (vertices of $G$) of $\pi^*$ are denoted
%by $(x,\jmath)$ where $x \in \Z$, $\jmath \in I$.
%Let 
%\[
%0 = \jmath_0 \prec \jmath_1 \prec \cdots \prec \jmath_s =M
%\]
%be the distinct values of the
%second components of the elements of $\pi^*$ in order of
%appearance in $\pi^*$.
%The sequence $(\jmath_0, \jmath_1, \ldots, \jmath_s)$ is not
%necessarily a path in $H(I, \preceq)$. Consider
%then another sequence $(\iota_0, \iota_1, \ldots, \iota_r)$
%which is a path in $H(I, \preceq)$ with $\iota_0=0$, $\iota_r=M$
%and 
%\[
%(\jmath_0, \jmath_1, \ldots, \jmath_s) \subset 
%(\iota_0, \iota_1, \ldots, \iota_r).
%\]

Let 
\[
0 = \iota_0 \prec \iota_1 \prec \cdots \prec \iota_s =M
\]
be the distinct values of the $I$-components of the 
elements of $\pi^*$ in order of appearance in $\pi^*$.
(The sequence $(\iota_0, \iota_1, \ldots, \iota_s)$ is not
necessarily a path in $H(I, \preceq)$.)
So for each $k=0,\ldots, s-1$, there are vertices 
$(x_k, \iota_k)$, $(y_k, \iota_{k+1})$ which are consecutive in the 
path $\pi^*$.
\fig{h}{9cm}{12cm}{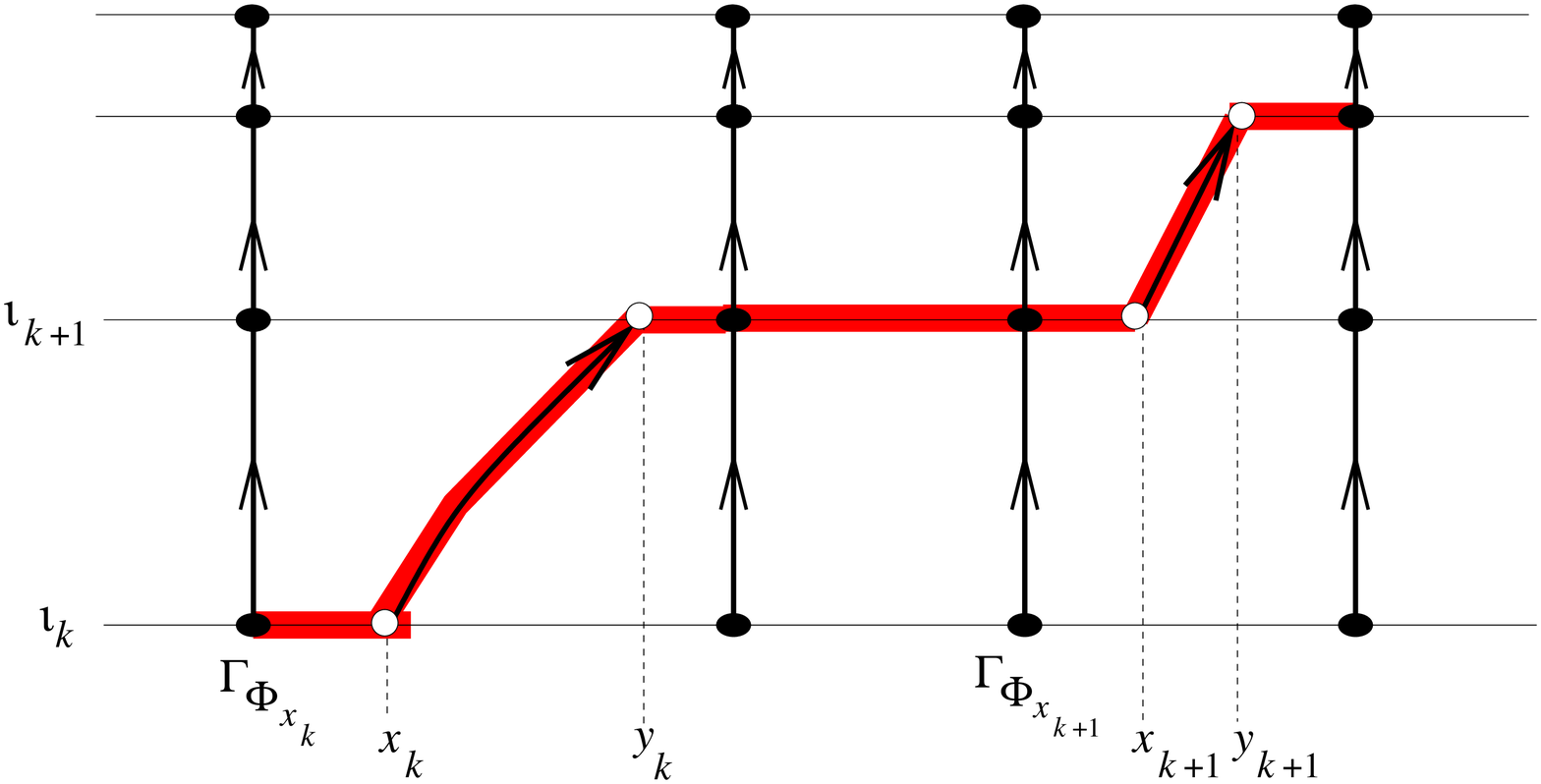}
{Construction used in obtaining the upper bound.}
Hence
\[
x_k \le y_k \le x_{k+1} , \quad \text{ for all } k=0,1,\ldots, s-1,
\]
where, by convention, we set $x_s=n$.
The point of $\SS$ prior to $x_k$ is $\Gamma_{\Phi_{x_k}}$ and,
since $\pi^*$ has maximum length, 
$(\Gamma_{\Phi_{x_k}}, \iota_k)$ is an element of $\pi^*$.
By the maximality of $\pi^*$ again, we have that $x_k$ and $y_k$ are
contained between two successive points of $\SS$ (otherwise we would
be able to strictly increase the length of the path).
Hence 
\begin{equation}
\label{ineqs}
\Gamma_{\Phi_{x_k}} \le x_k \le y_k \le \Gamma_{1+\Phi_{x_k}}
\le x_{k+1} , \quad \text{ for all } k=0,1,\ldots, s-1.
\end{equation}
We thus have
\[
L_n = |\pi^*| 
= L^{(\iota_0)}[0,\Gamma_1]+ L^{(\iota_0)}[\Gamma_1, \Gamma_{\Phi_{x_0}}]
+ \sum_{k=0}^{s-1} \bigg\{
L^{(\iota_k)}[\Gamma_{\Phi_{x_k}}, x_k] + 1 
	+ L^{(\iota_{k+1})}[y_k, \Gamma_{\Phi_{x_{k+1}}}] \bigg\}
+ L^{(\iota_s)}[\Gamma_{\Phi_n},n].  
\]
Due to \eqref{ineqs}, we have
\begin{align}
L^{(\iota_k)}[\Gamma_{\Phi_{x_k}}, x_k]
& \le L^{(\iota_k)}[\Gamma_{\Phi_{x_k}}, \Gamma_{1+\Phi_{x_k}}],
\label{nilk}
\\
L^{(\iota_{k+1})}[y_k, \Gamma_{\Phi_{x_{k+1}}}]
& \le L^{(\iota_{k+1})}[\Gamma_{\Phi_{x_k}}, \Gamma_{\Phi_{x_{k+1}}}],
\quad k=0, \ldots, s-1.
\label{good}
\end{align}
Moreover,
\begin{align}
L^{(\iota_0)}[0,\Gamma_1] & \le L^{(\iota_0)}[\Gamma_0,\Gamma_1]
\label{nil0}
\\
L^{(\iota_s)}[\Gamma_{\Phi_n},n] 
& \le L^{(\iota_s)}[\Gamma_{\Phi_n},\Gamma_{1+\Phi_n}].
\label{nils}
\end{align}
Each of the right-hand sides of \eqref{nilk}, \eqref{nil0} and
\eqref{nils} is bounded above by $\max_{0 \le j \le \Phi_n}
L^{(\iota)}[\Gamma_j, \Gamma_{1+j}]$.
If we then define
\[
\zeta_n := \sum_{\iota \in I} \max_{0 \le j \le \Phi_n}
L^{(\iota)}[\Gamma_j, \Gamma_{1+j}]
\]
and use \eqref{good}, we obtain
\[
L_n \le \zeta_n + M + 
\sum_{k=0}^{s-1} 
L^{(\iota_{k+1})}[\Gamma_{\Phi_{x_k}}, \Gamma_{\Phi_{x_{k+1}}}].
\]
Since for each sequence $0=\iota_0 \prec\iota_1 \prec \cdots \prec \iota_s=M$
of distinct ordered elements of $I$ we can find a path in the Hasse diagram
containing these elements, it follows easily that
\[
L_n \le \zeta_n + M + L_{\Phi_n}^*,
\]
which gives the upper bound. 
The upper bound is close to $L_n$
in the sense that the sequence the $\zeta_n$ are of order 1 in
distribution, i.e.\ that $(\zeta_n)$ is tight random sequence.
On the other hand, $nt = \Gamma_{\Phi_{[nt]}}-\Gamma_1 + o_{\mathrm d}(1)$.
It is thus clear that the weak limit of $\ell_n$ and that of 
\[
\ell^*_n(t) := 
\frac{L^*_{\Phi_{[nt]}}-C(\Gamma_{\Phi_{[nt]}}-\Gamma_1)}
{\kappa\sqrt{n}}, \quad t \ge 0,
\]
if it exists, will be identical.
Setting
\[
\ell^{*\!*}_n(t) := \frac{L^*_{[nt]}-C(\Gamma_{[nt]}-\Gamma_1)}{\kappa\sqrt{n}},
\quad \phi_n(t) := \frac{\Phi_{[nt]}}{n},
\]
we have 
\begin{equation}
\label{stars}
\ell^*_n(t)= \ell^{*\!*}_n(\phi_n(t)),
\end{equation}
and so the weak limit of $\ell_n^*$ is equal to that of $\ell_n^{*\!*}$
(if this exists) composed by the function $\{\lambda t\}$.

To show that the weak limit of $\ell_n^*$  exists and find it, define the
function $\psi : D[0,\infty)^I \to D[0,\infty)$ by
\begin{multline*}
\psi(\beta^{(i)}, i \in I)(t) :\max_{\iota}~
\sup_{\substack{0 \le t_0 \le t_1 \le \cdots \le t_{r}=t\\|\iota|=r}}
\bigg\{
\beta^{(\iota_0)}(t_{0})
+ \big[\beta^{(\iota_1)}(t_{1})-\beta^{(\iota_1)}(t_{0})\big]
+ \cdots
\\ \cdots
+ \big[\beta^{(\iota_{r})}(t_{r})-\beta^{(\iota_{r})}(t_{r-1})\big]
\bigg\}
\end{multline*}
where the maximum is taken over all paths $\iota$ 
from the minimum to the maximum element in the Hasse diagram
$H(I, \preceq)$.
The function $\psi$ is continuous (with respect to the topology
of uniform convergence).
Let
\[
s^{(i)}_n(t) := \frac{L^{(i)}[\Gamma_1, \Gamma_{[nt]}]-C(\Gamma_{[nt]}-\Gamma_1)}{\sigma \sqrt{n}}, \quad t \ge 0, \quad i \in I,
\]
where
\[
\sigma^2 := \var\{ L^{(i)}[\Gamma_1, \Gamma_2]-C(\Gamma_2-\Gamma_1)\}
\]
Since $L^{(i)}[\Gamma_j, \Gamma_{j+1}]$, $j \ge 1$, $i \in I$
are i.i.d.\ with common variance $\sigma^2$ we have (Theorem \ref{Lclt}) that
\begin{equation}
\label{BMvector}
\big( s_n^{(i)}, i \in I \big) \Rightarrow 
\big( B^{(i)}, i \in I \big)
\end{equation}
where $B^{(i)}$, $i \in I$ are i.i.d.\ standard Brownian motions.
Let
\[
\kappa := \lambda^{1/2} \sigma
\]
and observe that
%\footnote{An extra line or two of Algebra would be useful.}
\[
\ell_n^{*\!*}(t) = \lambda^{-1/2} \cdot \psi( s_n^{(i)}, i \in I )(t).
\]
By \eqref{BMvector} and the invariance principle,
\[
\ell_n^{*\!*} \Rightarrow \lambda^{-1/2} \cdot \psi(B^{(i)}, i \in I).
\]
By the relation \eqref{stars} and the remark following it,
we have
\[
\ell_n^* \Rightarrow  \psi(B^{(i)}, i \in I),
\]
and the right-hand side is equal in distribution to $Z$ (defined
by \eqref{Z1}-\eqref{Z2}). 
\qed

The remarks at the end of Section \ref{clt} also apply in the current 
case. We can easily conclude that $T_n$, the maximum length
of all paths from $(0,0)$ to $(n,M)$, has the same
asymptotics as $L_n$. In particular, Theorem \ref{ccclt}
holds if we replace $L_n$ by $T_n$. 

\section{Connection to last passage percolation}
Consider now the case
\[
I=\{0,1,\ldots, M\}
\]
with the usual ordering.
Assumption {\sf [D3]} can be substituted by
\\[5mm]
\hspace*{1cm}
\begin{minipage}{\textwidth}
{\sf [D3'] }~ 
For all $1 \le i \le M$ we have $r_{0,i-1,i} > 0$.
\end{minipage}
\\[5mm]
Let $G_M$ be the corresponding random directed cylinder graph, referred
to as slab graph here.
In particular, we can think of $G_M$ as the restriction of
a graph $G_{\infty}$ on the vertex set $\Z \times \Z_+$ where
two vertices $(x,i)$ and $(y,j)$, with $(x,i) \before (y,j)$, are connected
with probability $p_{y-x,j-i}$ that depends on the relative position
of the two vertices on the 2-dimensional lattice.

The problem here becomes that of a last passage percolation ,
%\footnote{Last
%passage percolation is a phenomenon and, by extension, a model. It is not
%a random variable as was written in the previous version.}, 
although
the model is not the standard nearest-neighbour one. Physically, we
can think of tunnels which run upwards (or
in directions southwest to northeast) and fluid moving in tunnels.
It takes one unit of time to cross a specific tunnel. We are interested
in the particle that starts from $(0,0)$ and reaches $(n,M)$ in the largest
possible time.
Since the Hasse diagram of
the set $\{0,1,\ldots, M\}$ with the natural ordering is the
linear graph with edges from $i-1$ to $i$, $1 \le i \le M$,
the limit process $Z$ is given by the simplified expression
\[
Z_t = \max_{0 \le t_0 \le \cdots \le t_M=t} 
\big\{B^{(0)}(t_0)+[B^{(1)}(t_1)-B^{(1)}(t_0)] + \cdots +
[B^{(M)}(t_M)-B^{(M)}(t_{M-1})]\big\}, \quad t \ge 0.
\]
The latter process is a Brownian last passage percolation process.
As was shown in  \cite{BAR2001, GTW2001, OCONYOR2002} it 
is a non-Gaussian process with marginal distribution
\[
Z_t \eqdist \sqrt{t} \cdot \lambda_M,
\]
for each $t \ge 0$, where $\lambda_M$ is the largest eigenvalue of a random
$(M+1)\times (M+1)$ Gaussian Unitary Ensemble (GUE) \cite{MEHTA04}. 

%GUE is defined in the space of Hermitian matrices as follows. 
%The matrix non-diagonal elements $\{\xi_{ij}+\sqrt{-1}~\eta_{ij}\}_{1\le i<j\le M}$  are 
%complex  random variables and diagonal elements $\{\xi_{ii}\}_{1\le i\le M} $ are real random variables.
%These random variables $\{\xi_{ij},\eta_{ij},\xi_{ii}\}_{1\le i<j\le M}$ are independent of each other and 
%have Gaussian distribition 
%with mean $0$  and variance
%$$
%Var (\xi_{ij})=Var (\eta_{ij})=1/2,\quad Var (\xi_{ii})=1.
%$$
%The Gaussian Unitary Ensemble has been extensively studied for the last 50 years, see \cite{MEHTA04}. 

Tracy and Widom \cite{TW93, TW94} showed that, as $M \to \infty$, the following
weak limit holds:
$$
M^{1/6}(\lambda_M-2\sqrt M)\Rightarrow F_{\text{TW}}, 
$$
with $F_{\text{TW}}$ being the Tracy-Widom distribution
whose hazard rate equals $\int_t^\infty q(x)^2 dx$, where $q(x)$
satisfies a Painlev\'e II equation; see \cite[eq.\ (3.1.7)]{AGZ}.
For an account on the universality of this distribution, see, e.g.,
\cite{DEI06}.
A number of interesting results have been proved relating this 
limiting distribution with certain stochastic models. 
These models include  longest increasing subsequence \cite{BDJ99}, 
last passage percolation, non-colliding particles, tandem queues 
\cite{BAR2001,GTW2001}, and random tilings \cite{JOH05}.
For the last passage percolation, in particular, this limit 
is known to appear in two cases. 
The first is the Brownian last passage percolation. 
The second is the last passage percolation model 
with exponential (or geometric) weights. 
In this model one puts independent and identically distributed 
exponential random variables 
in the vertices of $\Z_+^2$ and considers 
the maximum $L(M,N)$ of the sums of the weights 
over all directed paths from $(0,0)$ to $(M,N)$. 
It was shown in \cite{JOHANSSON00}
the random variable $L(N,N)$, properly normalized,
converges to the Tracy-Widom distribution as $N$ goes to infinity.
%What happens when the weights are not exponential 
%is essentially an open question. 
In \cite{BM2005}, more general weights were considered
and an analogous result for the random variable $L(N, N^a)$ (for
an appropriate $a$ depending on moment conditions) was obtained.
It is then natural to conjecture that a similar phenomenon occurs
in our slab graph too.

\section*{Acknowledgments}
We thank the Isaac Newton Institute for Mathematical Sciences for
providing the stimulating research atmosphere
where this research work was completed.
We are also grateful to Svante Janson and Graham Brightwell for pointing out
reference \cite{ABBJ94} to us.

\small
\vspace*{1cm}

\noindent
\begin{minipage}[t]{7cm}
\small \sc
Denis Denisov\\
School of Mathematics\\
Cardiff University\\
Cardiff CF24 4AG, UK\\
E-Mail: {\tt dennissov@googlemail.com}
\end{minipage}
\begin{minipage}[t]{6cm}
\small \sc
Serguei Foss \\
School of Mathematical Sciences\\
Heriot-Watt University\\
Edinburgh EH14 4AS, UK\\
E-mail: {\tt foss@ma.hw.ac.uk}
\end{minipage}
\\[5mm]
\begin{minipage}[t]{6cm}
\small \sc
Takis Konstantopoulos\\
School of Mathematical Sciences\\
Heriot-Watt University\\
Edinburgh EH14 4AS, UK\\
E-mail: {\tt takiskonst@gmail.com}
\end{minipage}

\end{document}